\documentclass[aap,MSNbibl,dvips]{arximspdf}
\usepackage{graphicx}
\doi{10.1214/12-AAP889} 
\volume{23}
\issue{5}
\pubyear{2013}
\firstpage{1879}
\lastpage{1912}

\makeatletter

\newcommand{\pr}{\mathbb{P}}
\newcommand{\E}{\mathbb{E}}

\newcommand{\reals}{\mathbb{R}}

\newtheorem{theorem}{Theorem}
\newtheorem{lemma}{Lemma}
\newtheorem{proposition}{Proposition}
\newtheorem{corollary}{Corollary}
\makeatother

\begin{document}
\begin{frontmatter}

\title{On the rate of convergence to stationarity of the M/M/N queue in
the Halfin--Whitt regime}
\runtitle{M/M/N queue in the Halfin--Whitt regime}

\begin{aug}
\author[A]{\fnms{David} \snm{Gamarnik}\thanksref{t1}\ead[label=e1]{gamarnik@mit.edu}}
\and
\author[B]{\fnms{David A.} \snm{Goldberg}\corref{}\thanksref{t2}\ead[label=e2]{dgoldberg9@isye.gatech.edu}}
\thankstext{t1}{Supported by NSF Grant CMMI-0726733.}
\thankstext{t2}{Supported by a Department of Defense NDSEG fellowship.}
\runauthor{D. Gamarnik and D.~A. Goldberg}
\affiliation{Massachusetts Institute of Technology and Georgia
Institute of Technology}
\address[A]{Operations Research Center and\\
\quad Sloan School of Management\\
Massachusetts Institute of Technology\\
100 Main Street, E62-563\\
Cambridge, Massachusetts 02139\\
USA\\
\printead{e1}} 
\address[B]{H. Milton Stewart School of\\
\quad Industrial and Systems Engineering\\
Georgia Institute of Technology\\
765 Ferst Drive, Room 437 Groseclose\\
Atlanta, Georgia 30332-0205\\
USA\\
\printead{e2}}
\end{aug}

\received{\smonth{3} \syear{2010}}
\revised{\smonth{7} \syear{2012}}

%
\begin{abstract}
We prove several results about the rate of convergence to stationarity,
that is, the spectral gap, for the $M/M/n$ queue in the Halfin--Whitt
regime. We identify the limiting rate of convergence to steady-state,
and discover an asymptotic phase transition that occurs w.r.t. this
rate. In particular, we demonstrate the existence of a constant $B^*
\approx1.85772$ s.t. when a certain excess parameter $B \in(0,B^*]$,
the error in the steady-state approximation converges exponentially
fast to zero at rate $\frac{B^2}{4}$. For $B > B^*$, the error in the
steady-state approximation converges exponentially fast to zero at a
different rate, which is the solution to an explicit equation given in
terms of special functions. This result may be interpreted as an
asymptotic version of a phase transition proven to occur for any fixed
$n$ by van Doorn [\textit{Stochastic Monotonicity and Queueing Applications of Birth-death
Processes} (1981) Springer].

We also prove explicit bounds on the distance to stationarity for the
$M/M/n$ queue in the Halfin--Whitt regime, when $B < B^*$. Our bounds
scale independently of $n$ in the Halfin--Whitt regime, and do not
follow from the weak-convergence theory.
\end{abstract}

%

\end{frontmatter}
%
\section{Introduction}\label{introsec}

Parallel server queueing systems can operate in a variety of regimes
that balance between efficiency
and quality of offered service. This is captured by the so-called
Halfin--Whitt (HW) heavy-traffic regime, which can be described as
critical w.r.t. the probability that an arriving job has to wait for
service. Namely, in this regime the stationary probability of wait is
bounded away from both zero and unity, as the number of servers grows.
Although studied originally by Pollaczek~\cite{Poll46} (see also
\cite
{JansenLRoots08}), Erlang~\cite{E48} and Jagerman~\cite{J74}, the
regime was formally introduced by Halfin and Whitt~\cite{HW81}, who
studied the $GI/M/n$ system for large $n$ when the traffic intensity
scales like $1 - Bn^{-{1}/{2}}$ for some strictly positive excess
parameter~$B$. They proved that, under minor technical assumptions on
the inter-arrival distribution, this sequence of $GI/M/n$ queueing
models has the following properties:
\begin{longlist}[(iii)]
\item[(i)] the steady-state probability that an arriving job has to wait for
service has a nontrivial limit;
\item[(ii)] the sequence of queueing processes, normalized by $n^{{1}/{2}}$, converges weakly to a nontrivial positive recurrent
diffusion, a.k.a. the HW diffusion;
\item[(iii)] the sequence of steady-state queue length distributions,
normalized by $n^{{1}/{2}}$, is tight and converges
distributionally to the mixture of a point mass at zero and an
exponential distribution.
\end{longlist}

Since the steady-state behavior of the $M/M/n$ queue in the HW regime
is quite simple~\cite{HW81}, while the transient dynamics are more
complicated~\cite{HW81}, it is common to use the steady-state
approximation to the transient distribution~\cite{GKM03}. Thus it is
important to understand the quality
of the steady-state approximation. The only work along these lines
seems to be the recent papers~\cite{KL08,KL10aa}, in which the
authors study the Laplace transform of the HW and related diffusions,
and prove several results analogous to our own for these diffusions.
The key difference is that in this paper we study the pre-limit
diffusion-scaled $M/M/n$ queue, not the limiting diffusion. We note
that the relevant transform functions were also studied in \cite
{AC68}, although in a different context. Also, similar questions were
studied for the associated sequence of fluid-scaled queues in~\cite{KR10}.

The question of how quickly the positive recurrent $M/M/n$ queue
approaches stationarity has a rich history in the queueing literature.
In~\cite{M55}, Morse derives an explicit solution for the transient
$M/M/1$ queue, and discusses implications for the exponential rate of
convergence to stationarity. Similar analyses are carried out in \cite
{C56} and~\cite{S60}. Around the same time, both Ledermann and Reuter
\cite{LR54}, and Karlin and McGregor (KM)~\cite{KM57}, worked out
powerful and elegant theories that could be used to give the transient
distributions for large classes of birth--death processes (BDP),
including the $M/M/n$ queue. The transient probabilities are expressed
as integrals against a spectral measure $\phi$, which is intimately
related to the eigenvalues of the generator of the BDP. KM devote an
entire paper~\cite{KM58} to the application of their theory to the
$M/M/n$ queue, in which they comment explicitly on the relationship
between the rate of convergence to stationarity and the support of
$\phi$. This relationship was later formalized in a series of papers by
other authors~\cite{C91,D02a}. Let $P(t)$ denote the matrix of
transient probabilities for the $M/M/n$ queue; that is, $P_{i,j}(t)$ is
the probability that there are $j$ jobs in system at time $t$, if there
are $i$ jobs in system at time 0. Let $A$ denote the generator matrix
associated with the $M/M/n$ queue, that is, $\frac{d}{dt} P(t) = A
\cdot P(t)$~\cite{EK05}. Recall that the spectral gap $\gamma$ of
a\vadjust{\goodbreak}
BDP is the absolute value of the supremum of the set of strictly
negative real eigenvalues of $A$ over an appropriate domain, and we
refer the reader to~\cite{C91} for details. Then we have the following
from the results of~\cite{C91}:
%
\begin{theorem}\label{background1}
For any positive recurrent $M/M/n$ or $M/M/\infty$ queue, $\gamma\in
(0, \infty)$. For all $i$ and $j$, $\lim_{t \rightarrow\infty} -
t^{-1} \log|P_{i,j}(t) - P_j(\infty)|$ exists, and is at least~$\gamma
$. For at least one pair of $i$ and $j$,
$\lim_{t \rightarrow\infty}
- t^{-1} \log|P_{i,j}(t) - P_j(\infty)| = \gamma$. Furthermore,
$\gamma= \inf\lbrace x\dvtx x > 0, \phi(x + \varepsilon) - \phi(x -
\varepsilon)
> 0 \mbox{ for all } \varepsilon> 0 \rbrace$.
\end{theorem}

We note that $\gamma$ is closely related to the singularities of the
Laplace transform of~$\phi$, and refer the reader to~\cite{KM58} for
details. It is well known that for the positive recurrent $M/M/1$ and
$M/M/\infty$ queues, $\gamma$ can be computed explicitly. In
particular, the following is proven in~\cite{KM58}:
%
\begin{theorem}\label{mm1infinity}
For the positive recurrent $M/M/1$ queue with arrival rate $\lambda$
and service rate $\mu$, $\gamma= (\lambda^{{1}/{2}}-\mu^{{1}/{2}})^2$, and
the spectral measure $\phi$ consists of a jump at zero, and an
absolutely continuous measure on $[(\lambda^{{1}/{2}}-\mu^{{1}/{2}})^2, (\lambda^{{1}/{2}}+\mu^{{1}/{2}})^2]$. For the
$M/M/\infty$ queue with arrival rate $\lambda$ and service rate $\mu$,
$\gamma= \mu$, and
the spectral measure $\phi$ consists of a countably infinite number of
jumps, with exactly one jump at every nonnegative integer multiple of
$\mu$.
\end{theorem}

Unfortunately, for the general positive recurrent $M/M/n$ queue, the
known characterizations for $\gamma$ involve computing the roots of
high-degree polynomials, which may be computationally difficult. This
arises from the fact that for the positive recurrent $M/M/n$ queue with
arrival rate $\lambda$ and service rate $\mu$, the spectral measure $\phi$ consists of three parts, as described in~\cite{KM58}. The first
part is a jump at zero, which
corresponds to the steady-state distribution. The second component is
an absolutely continuous measure on the
interval $ [  (\lambda^{{1}/{2}} - (n \mu)^{{1}/{2}}
)^2,  (\lambda^{{1}/{2}} + (n \mu)^{{1}/{2}} )^2
]$.
The third component consists of a set of at most $n$ (but possibly
zero) jumps, which all exist on $ ( 0,  ( \lambda^{{1}/{2}}
-
(n \mu)^{{1}/{2}}  )^2  )$.
The complexity of determining~$\gamma$ arises from the difficulty of
locating these jumps~\cite{D85}. In~\cite{KM58}, this set of jumps is
expressed in terms of the zeros of a certain polynomial equation.

Significant progress toward understanding these jumps was made in a
series of papers by van Doorn
\cite{D81,D84,D85,D87}. Van Doorn used the KM
representation and the theory of orthogonal polynomials to
give several alternate characterizations and bounds for the spectral
gap of a BDP, and applied these to the $M/M/n$ queue. He also
showed in~\cite{D81} that for each fixed $n$ there is a transition in
the nature of the spectral measure of the $M/M/n$ queue as one varies
the traffic intensity, proving the following theorem:

\begin{theorem}\label{fixedntransition}
For all $n \geq1$, there exists $\rho^*_n \in[0, 1)$ s.t. for
any $M/M/n$ queue with traffic intensity at least $\rho^*_n$, $\gamma=
( \lambda^{{1}/{2}} - (n \mu)^{{1}/{2}}  )^2$; and for
any $M/M/n$ queue with traffic intensity strictly less than $\rho^*_n$,
$\gamma<  ( \lambda^{{1}/{2}} - (n \mu)^{{1}/{2}}  )^2$.\vadjust{\goodbreak}
\end{theorem}

Unfortunately, all of the characterizations (including that of $\rho^*_n$) given by van Doorn are again stated in terms of the roots of
high-degree polynomials, and van Doorn
himself comments in~\cite{D85} that one is generally better off using
the approximations that he gives in the same paper. Van Doorn's work
was later extended by Kijima in~\cite{K92}, and similar results were
achieved by Zeifman using different techniques in~\cite{Z91}. It was
also shown in
\cite{Z91} that $\rho^*_n \leq(1-\frac{1}{n})^2$.

There are also some results in the literature for explicitly bounding
the distance to stationarity, as opposed to just identifying the
exponential rate of convergence. In~\cite{Z91}, Zeifman used tools
from the theory of differential equations to give explicit bounds on
the total variational distance between the transient and steady-state
distributions of a BDP, and explicitly examines the $M/M/n$ queue.
In~\cite{DZ09,DZP09}, van Doorn and Zeifman used the
techniques developed in~\cite{Z91} to derive explicit bounds on the
distance to stationarity for a different queueing model, and examined
how their bounds perform in a certain heavy-traffic regime (not HW). In
\cite{C98}, Chen developed very general bounds for the distance to
stationarity for Markov chains, and then applied these to BDP. However,
these bounds are generally not studied in the HW regime, and thus may
not scale desirably with $n$ in the HW regime. We note that the
complexity of bounding the distance to stationarity uniformly for a
sequence of BDP is related to the cutoff phenomenon for Markov chains~\cite{D96},
which has been studied in the context of queueing systems
\cite{FRT99}.

In this paper, we prove several results about the rate of convergence
to stationarity for the $M/M/n$ queue in the HW regime. We identify the
limiting rate of convergence to steady-state, that is, the spectral
gap, and discover an asymptotic phase transition that occurs w.r.t.
this rate. Specifically, let $\gamma_n$ denote the spectral gap
associated with the $M/M/n$ queue with arrival rate $n - B n^{{1}/{2}}$
and service rate equal to unity. Then we demonstrate the
existence of a constant $B^* \approx1.85772$ s.t. when the excess
parameter $B \in(0,B^*]$, $\lim_{n \rightarrow\infty} \gamma_n =
\frac
{B^2}{4}$. For $B > B^*$, $\lim_{n \rightarrow\infty} \gamma_n$
exists, and can be given as the solution to an explicit equation
involving special functions. This result may be interpreted as an
asymptotic version of the phase transition proven to occur for any
fixed $n$ by van Doorn in~\cite{D81}. Indeed, we prove that $\lim_{n
\rightarrow\infty} n^{{1}/{2}}(1 - \rho_n^*) = B^*$. It thus
follows from the results of~\cite{D81} (see Theorem \ref
{fixedntransition}) that $\gamma_n =  ( n^{{1}/{2}} - (n - B
n^{{1}/{2}})^{{1}/{2}}  )^2$ for $B < B^*$ and all
sufficiently large~$n$. Observing that $\lim_{n \rightarrow\infty}
( n^{{1}/{2}} - (n - B n^{{1}/{2}})^{{1}/{2}}  )^2 =
\frac{B^2}{4}$ links our results to those of van Doorn for the case $B
< B^*$, and a similar connection exists for the case $B \geq B^*$.

We also prove explicit bounds on the distance to stationarity
for the $M/M/n$ queue in the HW regime, when $B < B^*$. Our bounds
scale independently of $n$ in the HW regime, and do not follow from the
weak-convergence theory.

\subsection{Outline of the paper}
The rest of the paper proceeds as follows. In Section~\ref{1mainsec},
we state our main results, and outline our proof technique.
In Section~\ref{1techprelim}, we prove a new characterization for the
spectral gap of the $M/M/n$ queue. In Sections~\ref
{1asssection}--\ref
{1zetazerosec}, we study the asymptotic properties of this
characterization. In Section~\ref{1gammansec}, we compute the limiting
spectral gap of the $M/M/n$ queue in the HW regime, and prove that a
phase transition occurs. In Section~\ref{explicitsec}, we prove our
explicit bounds on the distance to stationarity. In Section \ref
{comparesec}, we compare our explicit bounds to other bounds from the
literature. In Section~\ref{2conc} we summarize our main results and
present ideas for future research. We include a technical appendix in
Section~\ref{2appsec2}.

\section{Main results}\label{1mainsec}

\subsection{Definitions and notation}

Let ${\mathcal Q}^n$ denote the $M/M/n$ queue with arrival rate
$\lambda_n \stackrel{\Delta}{=} n - B n^{{1}/{2}}$ and service rate $\mu
\stackrel{\Delta}{=} 1$, where we assume throughout that $n$ is
sufficiently large to ensure that $\lambda_n > 0$, and $n > \lambda_n +
1$. Let $Q^n(t)$ denote the number in system, that is, the number of
jobs in service plus the number of jobs waiting in queue, at time $t$;
$Q^n(\infty)$ denote the corresponding steady-state r.v.; and $\gamma_n$ denote the spectral gap of the associated Markov chain. We define
$P^n_{i,j}(t) \stackrel{\Delta}{=} \pr ( Q^n(t) = j | Q^n(0) = i
)$, $P^n_j(\infty) \stackrel{\Delta}{=} \pr ( Q^n(\infty)
= j
)$, $P^n_{i, \leq j}(t) \stackrel{\Delta}{=} \sum_{k=0}^j
P^n_{i,k}(t)$ and
$P^n_{\leq j}(\infty) \stackrel{\Delta}{=}  \sum_{k=0}^j
P^n_{k}(\infty
)$. For a function $f$, we let $Z(f)  ( Z^+(f)  )$ denote the
infimum of the set of (strictly positive) real zeros of $f$, and set $
Z(f)  ( Z^+(f)  ) = \infty$ if $f$ has no (strictly positive)
real zeros. All logarithms will be base $e$. Unless otherwise stated,
all functions are defined only over $\reals$. All empty products are
assumed to be equal to unity, and all empty summations are assumed to
be equal to zero. Also, for an event $\lbrace E \rbrace$, we let
$I
(\lbrace E \rbrace )$ denote the corresponding indicator function.

\subsection{The parabolic cylinder functions}
We now briefly review the two-parameter function commonly referred to
as the parabolic cylinder function $D_x(z)$, since we will need
these functions for the statement (and proof) of our main results. For
excellent references on these functions, see
\cite{G80} Sections~8.31 and~9.24,~\cite{B69}~Sections
3.3--3.5 and~\cite{E53} Chapter 8.
Let $\Gamma$ denote the Gamma function (see~\cite{H59}, Chapter 8.8).
It is stated in~\cite{B69} that $x,z \in\reals$ implies $D_x(z) \in
\reals$, and
%
\begin{equation}
\label{gonefishin2}\quad  D_x(z) = %
\cases{\displaystyle\biggl(
\frac{2}{\pi}\biggr)^{{1}/{2}} \exp\biggl( \frac{z^2}{4} \biggr) \int
_{0}^{\infty} \exp\biggl( -\frac{y^2}{2} \biggr)
\cos\biggl( \frac{\pi}{2} x - z y \biggr) y^x \,dy,\vspace*{2pt}\cr
\hspace*{227pt}\mbox{if } x
\geq0,\vspace*{2pt}
\cr
\displaystyle\frac{ \exp( - {z^2}/{4} ) }{\Gamma(-x)} \int_{0}^{\infty}
\exp\biggl( -\frac{y^2}{2} - z y \biggr) y^{ -(x+1) } \,dy, \qquad  \mbox{if } x <
0.} %
\end{equation}
$D_x(z)$ takes on a simpler form for integral $x$. In
particular, it is stated in~\cite{G80} that for $z \in\reals$,
%
\begin{eqnarray}
\label{parainteger} D_{-1}(z) &=& 2^{{1}/{2}} \exp\biggl(
\frac{z^2}{4} \biggr) \int_{ 2^{-{1}/{2}} z }^{\infty} \exp\bigl(
-y^2 \bigr)\,dy,
\nonumber
\\[-8pt]
\\[-8pt]
\nonumber
  D_0(z) &=& \exp\biggl( -\frac{z^2}{4}
\biggr)\quad \mbox{and} \quad D_1(z) = z \exp\biggl( -\frac{z^2}{4}
\biggr).
\end{eqnarray}
 Note that since $\Gamma(-x) \in(0, \infty)$ for $x < 0$, (\ref{gonefishin2}) and (\ref{parainteger}) imply that
$D_x(z) > 0$ for $z \in\reals$ and $x \leq0$.

The parabolic cylinder functions arise in several contexts
associated with the limits of queueing models, such as the
Ornstein--Uhlenbeck limit of the appropriately scaled infinite-server
queue~\cite{I65} and various limits associated with the Erlang loss
model~\cite{XK93}. We note that the parabolic cylinder functions have
been studied as the limits of certain polynomials under the HW scaling,
using tools from the theory of differential equations \cite
{D07,D08,D01,BW94}.
\subsection{Main results}
\begin{figure}

\includegraphics{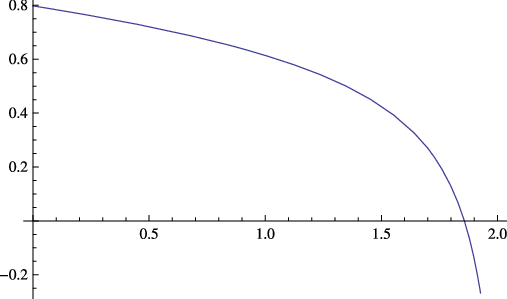}

\caption{Plot of $\zeta$.}
\label{figzetab}
\end{figure}

We now state our main results. We begin by identifying the limiting
rate of convergence to steady-state, that is, the limiting spectral
gap, for the $M/M/n$ queue in the HW regime; and prove that a phase
transition occurs w.r.t. this limiting rate. We define
\[
\upsilon(x,y) \stackrel{\Delta} {=} %
\cases{\displaystyle \frac{D_x(y) }{ D_{x-1}(y) }, &\quad  $
\mbox{if } D_{x-1}(y) \not= 0,$ \vspace*{2pt}
\cr
\infty,&\quad $
\mbox{otherwise.}$} %
\]
Also, let $\varphi(B) \stackrel{\Delta}{=} \upsilon(\frac{B^2}{4},
-B)$, $\zeta(B) \stackrel{\Delta}{=} \varphi(B) + \frac{B}{2}$ and
\[
\Psi_{\infty}(x) \stackrel{\Delta} {=} %
\cases{\displaystyle \upsilon(x,-B) +
\frac{1}{2} \bigl(B + \bigl(B^2 - 4x\bigr)^{{1}/{2}}
\bigr), &\quad $ \mbox{if }\displaystyle  x \leq\frac{B^2}{4},$ \vspace*{2pt}
\cr
\infty,&\quad  $
\mbox{otherwise.}$} %
\]
Note that $\zeta(B) = \Psi_{\infty}(\frac{B^2}{4})$. We include a plot
of $\zeta$ in Figure~\ref{figzetab}.
%

\begin{figure}[b]

\includegraphics{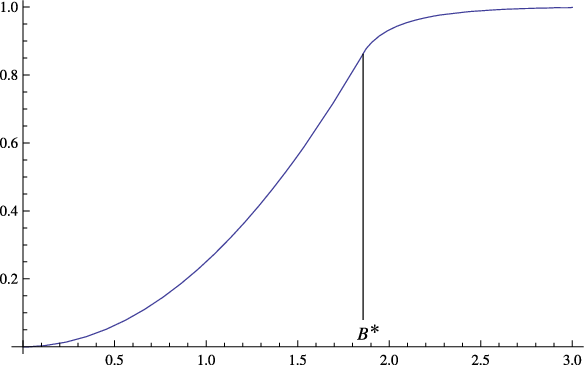}

\caption{Plot of $\gamma$.}\label{figgammab}
\end{figure}

Let $B^* \stackrel{\Delta}{=} Z^+( \zeta).$ Then:
%
\begin{proposition}\label{Bstarprops}
$B^* \approx1.85772$ and $Z^+( \Psi_{\infty} ) \in (0, \min(1,
\frac{B^2}{4})  )$ for \mbox{$B > B^*$}.
\end{proposition}
 Our main result is:
%
\begin{theorem}\label{asymptoticphasetransition}
The limit $\gamma(B) \stackrel{\Delta}{=} \lim_{n \rightarrow
\infty}
\gamma_{n}$ exists for all $B > 0$.
For $0 < B \leq B^*$,
$\gamma(B) = \frac{B^2}{4}.$
For $B \geq B^*$, $\gamma(B) = Z^+( \Psi_{\infty} )$.
\end{theorem}
We include a plot of $\gamma$ in Figure~\ref{figgammab}.\vadjust{\goodbreak}
%

Due to the nonlinear manner in which the steady-state probability of
wait scales in the HW regime, the case $0 <
B < B^*$ actually encompasses most scenarios of practical interest.
Indeed, it is proven in~\cite{HW81} that
the limit of the steady-state probability of wait equals
\[
\biggl(1 + B \exp\biggl(\frac{1}{2}B^2\biggr)\int
_{-\infty}^B \exp\biggl(-\frac{1}{2}z^2
\biggr) \,dz \biggr)^{-1}.
\]
As this limit is monotone in $B$, the case $0 < B < B^*$ includes all
scenarios for which the steady-state probability of wait is at least $0.04$.

We note that the results of~\cite{KL08} show that $\gamma
(B)$ is also the spectral gap of the HW diffusion, demonstrating an
interchange of limits for the $M/M/n$ queue in the HW regime.\vadjust{\goodbreak} Namely,
the limit of the sequence of spectral gaps equals the spectral gap of
the corresponding weak limit. Interestingly, neither result implies the
other, and it is an open challenge to understand this interchange more
generally.

The following corollary may be interpreted as an asymptotic version of
Theorem~\ref{fixedntransition}.
%
\begin{corollary}\label{asymptoticphasetransitioncor}
The $\rho^*_n$ parameter of Theorem~\ref{fixedntransition} satisfies
\[
\lim_{n \rightarrow\infty} n^{{1}/{2}} \bigl(1 - \rho^*_n\bigr) = B^*.
\]
\end{corollary}
 We now give an interpretation of Theorem \ref
{asymptoticphasetransition} and Corollary \ref
{asymptoticphasetransitioncor}. The $M/M/n$ queue behaves like an
$M/M/1$ queue when all servers are busy, and an $M/M/\infty$ queue when
at least one server is idle.
The phase transition of Theorem~\ref{asymptoticphasetransition}
formalizes this relationship in a new way.
For $0 < B < B^*$, the KM
spectral measure of the $M/M/n$ queue in the HW regime has no jumps
away from the origin, and has spectral gap equal to $( \lambda_n^{{1}/{2}} - n^{{1}/{2}} )^2$,
two properties shared by the associated
$M/M/1$ queue; see Theorem~\ref{mm1infinity}. For $B > B^*$, the KM
spectral measure has at
least one jump away from the origin, like the associated $M/M/\infty$
queue (whose spectral measure has only jumps and spectral gap
equal to unity; see Theorem~\ref{mm1infinity}). Another interpretation
is that the $M/M/n$ queue cannot approach stationarity faster than either
component system would on its own.

We now state our explicit bounds on the distance to stationarity
for the case $B < B^*$.
%
\begin{theorem}\label{explicitboundstheorem}
Given $B \in(0,B^*)$ and $a_1,a_2 \in\reals$, let $a = \max
(|a_1|,|a_2|,B)$. Then there exists $N_{B,a_1,a_2} < \infty$, depending
only on $B,a_1$ and $a_2$, s.t. for all $n \geq N_{B,a_1,a_2}$ and $t
\geq1$,
%
\begin{eqnarray}
\label{slayme11}&& \bigl|n^{{1}/{2}} P^n_{\lceil n + a_1 n^{{1}/{2}} \rceil,
\lceil n + a_2 n^{{1}/{2}} \rceil}(t) -
n^{{1}/{2}} P^n_{\lceil
n + a_2 n^{{1}/{2}}
\rceil}(\infty) \bigr|
\nonumber
\\[-8pt]
\\[-8pt]
\nonumber
&&\qquad\leq t^{-{1}/{2}}
\exp \biggl( 30 \bigl( a^2 + 1 \bigr) - \frac{B^2}{4} t \biggr)
\end{eqnarray}
and
%
\begin{eqnarray}
\label{slayme22} &&\bigl|P^n_{\lceil n + a_1 n^{{1}/{2}} \rceil, \leq\lceil n + a_2
n^{{1}/{2}} \rceil}(t) - P^n_{ \leq\lceil n + a_2 n^{{1}/{2}}
\rceil}(
\infty) \bigr|
\nonumber
\\[-8pt]
\\[-8pt]
\nonumber
&&\qquad \leq B^{-1} t^{-{1}/{2}}\exp \biggl( 30 \bigl(
a^2 + 1 \bigr) - \frac
{B^2}{4} t \biggr).
\end{eqnarray}
\end{theorem}

Note that Theorem~\ref{explicitboundstheorem} provides a bound for any
sufficiently large fixed $n$ and all times $t$ greater than unity,
which is independent of $n$, and converges to zero as $t \rightarrow
\infty$. Interestingly, such uniform bounds do not follow directly from
the weak-convergence theory, since the standard framework of weak\vadjust{\goodbreak}
convergence requires that one first fix a finite time interval of
interest, and then let $n \rightarrow\infty$, in that order.

It follows from the weak-convergence theory that our explicit
bounds yield corresponding bounds for the distance to stationarity of
the HW diffusion. Furthermore, in light of Theorem \ref
{asymptoticphasetransition}, the exponent $\frac{B^2}{4}$ appearing in
our bounds is the best possible. Although we were able to derive
partial results for the case $B \geq B^*$,
the derived bounds were considerably more complicated than those of
Theorem~\ref{explicitboundstheorem}, and we leave it as an open
question to derive simple explicit bounds for the case $B \geq B^*$. We
note that the results of~\cite{KL08} suggest that the exponential
dependence on $a^2$, and inverse dependence on $t^{{1}/{2}}$, of
the prefactor appearing in Theorem~\ref{explicitboundstheorem} may not
be tight, and it seems likely that a more refined analysis would yield
sharper bounds.

\subsection{Outline of proof}
We now present an outline of the proof of our main results. To prove
Theorem~\ref{asymptoticphasetransition} and Corollary \ref
{asymptoticphasetransitioncor}, we give a new characterization for the
spectral gap $\gamma_n$, and then study its asymptotics in the HW
regime. More precisely, in Section~\ref{1techprelim}, we prove a new
characterization for the spectral gap $\gamma_n$, in terms of a certain
function $\Psi_n$ which we define. We express $\gamma_n$ in terms of
three quantities: $(n^{{1}/{2}}-\lambda_n^{{1}/{2}})^2$,
$Z^+(\Psi_n)$ and the sign of $\Psi_n ((n^{{1}/{2}}-\lambda_n^{{1}/{2}})^2 )$.
In Section~\ref{1asssection}, we prove that
in the HW regime, $\Psi_n$ converges to $\Psi_{\infty}$ and $\Psi_n
((n^{{1}/{2}}-\lambda_n^{{1}/{2}})^2 )$ converges to
$\zeta
(B)$. In Section~\ref{1zplussec}, we prove that in the HW regime,
$Z^+(\Psi_n)$ converges to $Z^+(\Psi_{\infty})$.
In Section~\ref{1zetazerosec}, we characterize the sign of $\zeta(B)$.
In Section~\ref{1gammansec}, we combine the above results to prove
Theorem~\ref{asymptoticphasetransition} and Corollary \ref
{asymptoticphasetransitioncor}. To prove Theorem \ref
{explicitboundstheorem}, we use induction arguments to bound certain
polynomials which appear in the KM representation for the transient
$M/M/n$ queue.

\section{\texorpdfstring{Characterization for $\gamma_n$}
{Characterization for gamma n}}\label{1techprelim}
In this section we give a new characterization for~$\gamma_n$. We begin
by associating several functions to the $M/M/n$ queue, as in \cite
{K92} and~\cite{D81}. For $0 \leq k \leq n$, let $f_{n,k}(x)
\stackrel
{\Delta}{=} \sum_{j=0}^k {k \choose j} \lambda_n^j \prod_{i=1}^{k-j} (i-x)$;
\[
z_{n,k}(x) \stackrel{\Delta} {=} %
\cases{ \displaystyle\frac{f_{n,k}(x)}{f_{n,k-1}(x)},
&\quad $\mbox{if } f_{n,k-1}(x) \not= 0,$ \vspace*{2pt}
\cr
\infty,&\quad $
\mbox{otherwise,}$} %
\]
and $z_n(x) \stackrel{\Delta}{=} z_{n,n}(x).$
We also define
\[
a_n(x) \stackrel{\Delta} {=} %
\cases{\displaystyle \tfrac{1}{2}
\bigl( \lambda_n + n - x - \bigl( \bigl(n^{{1}/{2}} -
\lambda_n^{{1}/{2}}\bigr)^2 - x \bigr)^{{1}/{2}}
\bigl( \bigl(n^{{1}/{2}} + \lambda_n^{{1}/{2}}
\bigr)^2 - x \bigr)^{{1}/{2}} \bigr), \vspace*{2pt}\cr
\hspace*{38pt}\mbox{if } x \leq
\bigl(n^{{1}/{2}} - \lambda_n^{{1}/{2}}\bigr)^2,
\vspace*{2pt}
\cr
\infty,\qquad \mbox{otherwise} } %
\]
and
\[
\Psi_{n}(x) \stackrel{\Delta} {=} %
\cases{
z_n(x) - a_n(x), &\quad $\mbox{if } z_n(x)\not=
\infty \mbox{ or } a_n(x)\not= \infty,$ \vspace*{2pt}
\cr
\infty,& \quad $
\mbox{otherwise.}$}\vadjust{\goodbreak} %
\]
We now cite some properties of $f_{n,n-1}, z_{n,k}$ and $\Psi_n$, as
stated in~\cite{K92}, for use in later proofs.
%
\begin{lemma} \label{kijima2}
\textup{(i)} $f_{n,n-1}$ is strictly positive on $ (-\infty, 1 ]$.\vspace*{-6pt}
\begin{longlist}[(iii)]
\item[(ii)] For $k \leq n$, $z_{n,k}$ is strictly positive, continuous and
strictly decreasing on $( -\infty, 1 ]$.
\item[(iii)]$\Psi_n$ is continuous and strictly decreasing on $ (
-\infty,
\min ( (n^{{1}/{2}} - \lambda_n^{{1}/{2}})^2, 1
) ]$.
\end{longlist}
\end{lemma}
We now prove the main result of this section, a new characterization
for $\gamma_n$, in particular.
%
\begin{proposition}\label{crossatleastonce}
\textup{(i)} If $(n^{{1}/{2}} - \lambda_n^{{1}/{2}})^2 < 1$ and $\Psi_n
( (n^{{1}/{2}} - \lambda_n^{{1}/{2}})^2  ) < 0$, then
$\gamma_n = Z^+( \Psi_n )$.\vspace*{-6pt}

\begin{longlist}[(iii)]
\item[(ii)]
If $(n^{{1}/{2}} - \lambda_n^{{1}/{2}})^2 < 1$ and $\Psi_n
( (n^{{1}/{2}} - \lambda_n^{{1}/{2}})^2  ) \geq0$, then
$\gamma_n = (n^{{1}/{2}} - \lambda_n^{{1}/{2}})^2$.
\item[(iii)]
If $(n^{{1}/{2}} - \lambda_n^{{1}/{2}})^2 \geq1$, then $Z^+(
\Psi_n ) \in ( 0, 1  )$, and $\gamma_n = Z^+( \Psi_n
)$.
\end{longlist}
\end{proposition}
The proof of Proposition~\ref{crossatleastonce} relies on the
following known characterization for $\gamma_n$.
Let $\sigma_n(x) \stackrel{\Delta}{=} f_{n,n}(x) - (\lambda_n
n)^{{1}/{2}} f_{n,n-1}(x)$, and
$\psi_n(x) \stackrel{\Delta}{=} f_{n,n}(x) - a_n(x) f_{n,n-1}(x)$. Then
the following is proven in~\cite{K92}:
%
\begin{theorem} \label{kijimamainresult}
If $Z( \sigma_n ) \geq(n^{{1}/{2}} - \lambda_n^{{1}/{2}})^2$,
then $\gamma_n = (n^{{1}/{2}} - \lambda_n^{{1}/{2}})^2$. If
$Z( \sigma_n ) < (n^{{1}/{2}} - \lambda_n^{{1}/{2}})^2$, then
$\gamma_n = Z( \psi_n )$.
\end{theorem}
With Theorem~\ref{kijimamainresult} in hand, we now complete the proof
of Proposition~\ref{crossatleastonce}.
\begin{pf*}{Proof of Proposition~\ref{crossatleastonce}}
We begin by studying the sign of $\Psi_n(0),\break\Psi_n(1), \sigma_n(0)$ and
$\sigma_n(1)$. Note that
\begin{eqnarray*}
\Psi_n(0) &=& \frac{ \sum_{k=0}^n {n \choose k} \lambda_n^k (n-k)!
} {
\sum_{k=0}^{n-1} {n-1 \choose k} \lambda_n^k (n-1-k)! } - \frac
{1}{2} \bigl(
\lambda_n + n - \bigl((\lambda_n + n)^2 - 4
\lambda_n n \bigr)^{{1}/{2}} \bigr)
\\
&=& n \frac{ \sum_{k=0}^n { \lambda_n^k }/{k!} }{ \sum_{k=0}^{n-1} { \lambda_n^k }/{k!} } - \lambda_n > 0.
\end{eqnarray*}
If $( n^{{1}/{2}} - \lambda_n^{{1}/{2}} )^2 \geq1$, then
\begin{eqnarray*}
\Psi_n(1) &=& \frac{ \sum_{k=0}^n {n \choose k} \lambda_n^k \prod_{i=1}^{n-k} (i-1) } { \sum_{k=0}^{n-1} {n-1 \choose k}
\lambda_n^k \prod_{i=1}^{n-1-k} (i-1) } \\
&&{}- \frac{1}{2} \bigl(
\lambda_n + n - 1 - \bigl(( \lambda_n + n - 1
)^2 - 4 \lambda_n n \bigr)^{{1}/{2}} \bigr)
\\
&=& \frac{ \lambda_n^n }{ \lambda_n^{n-1} } - \frac{1}{2} \bigl( \lambda_n + n - 1
- \bigl(( \lambda_n + n - 1 )^2 - 4 \lambda_n
n \bigr)^{{1}/{2}} \bigr)
\\
&=& \frac{1}{2} \bigl( \lambda_n - n + 1 + \bigl( (
\lambda_n - n + 1)^2 - 4 \lambda_n
\bigr)^{{1}/{2}} \bigr) \leq 0.
\end{eqnarray*}
Similarly,
\begin{eqnarray*}
\sigma_n(0) &=& \sum_{k=0}^n
\pmatrix{n \cr k} \lambda_n^k (n-k)! - (\lambda_n
n)^{{1}/{2}} \sum_{k=0}^{n-1} \pmatrix{n-1
\cr k} \lambda_n^k (n-1-k)!
\\
&=& (n-1)! \Biggl( n \sum_{k=0}^n
\frac{ \lambda_n^k }{k!} - (\lambda_n n)^{{1}/{2}} \sum
_{k=0}^{n-1} \frac{ \lambda_n^k }{k!} \Biggr)\\
& \geq& (n-1)!
\sum_{k=0}^n \frac{ \lambda_n^k }{k!} \bigl( n -
(\lambda_n n)^{{1}/{2}} \bigr) > 0
\end{eqnarray*}
and
\begin{eqnarray*}
\sigma_n(1) &=& \sum_{k=0}^n
\pmatrix{n \cr k} \lambda_n^k \prod
_{i=1}^{n-k} (i-1) - (\lambda_n
n)^{{1}/{2}} \sum_{k=0}^{n-1} \pmatrix{n-1
\cr k} \lambda_n^k \prod_{i=1}^{n-1-k}
(i-1)
\\
&=& \lambda_n^n - (\lambda_n
n)^{{1}/{2}} \lambda_n^{n-1} < 0.
\end{eqnarray*}
We first prove assertion (i). Note that if $f_{n,n-1}(x) \not= 0$, then
$z_n(x) - (\lambda_n n)^{{1}/{2}} = \frac{\sigma_n(x)}{f_{n,n-1}(x)}$. Thus Lemma~\ref{kijima2}(i) implies
that $\sigma_n$ is the same sign as $z_n - (\lambda_n n)^{{1}/{2}}$
on $ (-\infty, (n^{{1}/{2}} - \lambda_n^{{1}/{2}})^2 ]$.
Recalling that $\sigma_n(0) > 0$, it follows from the
continuity/monotonicity of $z_n$ [guaranteed by Lemma \ref
{kijima2}(ii)] and the intermediate value theorem that
$\sigma_n$ has a zero on $ (-\infty, (n^{{1}/{2}} - \lambda_n^{{1}/{2}})^2  )$ if and only if $z_n  ( (n^{{1}/{2}} -
\lambda_n^{{1}/{2}})^2  ) - (\lambda_n n)^{{1}/{2}} <
0$. Since
$a_n  ( ( n^{{1}/{2}} - \lambda_n^{{1}/{2}} )^2  ) =
(\lambda_n n)^{{1}/{2}}$, we conclude that
$Z( \sigma_n ) < (n^{{1}/{2}} - \lambda_n^{{1}/{2}})^2$ iff
$\Psi_n  ( (n^{{1}/{2}} - \lambda_n^{{1}/{2}})^2  )
< 0$.
Thus $Z( \sigma_n ) < (n^{{1}/{2}} - \lambda_n^{{1}/{2}})^2$,
since by assumption $\Psi_n  ( (n^{{1}/{2}} - \lambda_n^{{1}/{2}})^2  ) < 0$, and $\gamma_n = Z(\psi_n)$ by Theorem \ref
{kijimamainresult}. Noting that $\Psi_n = \frac{ \psi_n }{ f_{n,n-1} }$
on $ (-\infty, (n^{{1}/{2}} - \lambda_n^{{1}/{2}})^2 ]$,
this further implies that $\gamma_n = Z( \Psi_n )$. That $\gamma_n =
Z^+(\Psi_n)$ then follows from the fact that $\Psi_n(0) > 0$, and the
continuity/monotonicity of $\Psi_n$ guaranteed by Lemma \ref
{kijima2}(iii). This completes the proof of assertion (i). The
proof of assertion (ii) follows from a similar argument, and we omit
the details.

We now prove assertion (iii). Since $\sigma_n$ is a
polynomial s.t. $\sigma_n(0) > 0$ and $\sigma_n(1) < 0$, we have that
$Z(\sigma_n) < 1 \leq(n^{{1}/{2}} - \lambda_n^{{1}/{2}})^2$.
Thus Theorem~\ref{kijimamainresult} implies that $\gamma_n = Z(\psi_n)$. As in the proof of assertion (i), it follows that $\gamma_n =
Z(\Psi_n)$. Since $\Psi_n(0) > 0$ and $\Psi_n(1) < 0$, the
continuity/monotonicity of $\Psi_n$ guaranteed by Lemma \ref
{kijima2}(iii) further ensures that $\gamma_n = Z^+(\Psi_n)
\in
(0,1)$, completing the proof.
\end{pf*}

\section{\texorpdfstring{Asymptotic analysis of $\Psi_n$}
{Asymptotic analysis of Psi n}}\label{1asssection}
In this section we derive the asymptotics of $\Psi_n$ in the HW regime.
In particular, we prove that:
%
\begin{theorem}\label{psiasymptotics}
For $B > 0$ and $x \in(0,1) \cap(0, \frac{B^2}{4}]$,
$\lim_{n \rightarrow\infty} \lambda_n^{-{1}/{2}} \Psi_n(x)
= \Psi_{\infty}(x).$
\end{theorem}
 We also prove that:
%
\begin{corollary}\label{asscorrect}
For $B \in(0,2)$, $\lim_{n \rightarrow\infty} \lambda_n^{-{1}/{2}} \Psi_n  ( (n^{{1}/{2}} - \lambda_n^{{1}/{2}})^2
) = \zeta(B)$.
\end{corollary}
 We proceed by separately analyzing the asymptotics of
$\lambda_n^{-{1}/{2}}  ( a_n - \lambda_n  )$ and $\lambda_n^{-{1}/{2}}  ( z_n - \lambda_n  )$, beginning with $a_n$. Let
\[
a_{\infty}(x) \stackrel{\Delta} {=} %
\cases{ \displaystyle\frac{1}{2}
\bigl(B - \bigl(B^2 - 4x\bigr)^{{1}/{2}} \bigr), & \quad$\mbox{if } x
\leq\displaystyle\frac{B^2}{4},$ \vspace*{2pt}
\cr
\infty,& \quad $\mbox{otherwise.}$} %
\]
Then:
%
\begin{lemma}\label{aasymptotics}
For $x \in[0,\frac{B^2}{4}]$, $\lim_{n \rightarrow\infty} \lambda_n^{-{1}/{2}}  ( a_n(x) - \lambda_n  ) = a_{\infty}(x)$.
\end{lemma}
\begin{pf}
Note that
\begin{eqnarray*}
&&\lambda_n^{-{1}/{2}} \bigl( a_n(x) -
\lambda_n \bigr)\\
&&\qquad = \bigl( Bn^{{1}/{2}} - x - \bigl(
\bigl(n^{{1}/{2}} + \lambda_n^{{1}/{2}}\bigr)^2
- x \bigr)^{{1}/ {2}} \bigl( \bigl(n^{{1}/{2}} -
\lambda_n^{{1}/{2}}\bigr)^2 - x \bigr)^{{1}/{2}} \bigr) \bigl(2 \lambda_n^{{1}/{2}}\bigr)^{-1}.
\end{eqnarray*}
The lemma then follows from the fact that $\lim_{n \rightarrow\infty}
( B n^{{1}/{2}} - x )(2 \lambda_n^{{1}/{2}})^{-1} = \frac{B}{2},
\lim_{n \rightarrow\infty}  ( ( n^{{1}/{2}} + \lambda_n^{{1}/{2}} )^2 - x  )^{ {1}/{2} }  ( 2 \lambda_n^{
{1}/{2} }
)^{-1} = 1$ and $\lim_{n \rightarrow\infty} ( n^{{1}/{2}} -
\lambda_n^{{1}/{2}} ) = \frac{B}{2}$.
\end{pf}
 We now analyze the asymptotics of $z_n$, and begin by proving
some necessary bounds.
Let us fix some $x \in(0,1)$ and integer $T \geq3$, and define
\begin{eqnarray*}
R_{1,n} &\stackrel{\Delta} {=} &\lambda_n^{{(x-1)}/{2}} \sum
_{k=0}^{n-(T+1)} (n - k)^{1-x} \exp(-
\lambda_n)\frac{\lambda_n^k}{k!}, \\
 R_{2,n} & \stackrel{\Delta}
{=} &\lambda_n^{{(x-2)}/{2}} \sum_{k=0}^{ \lceil n - T^{-1} n^{{1}/{2}} \rceil}
k (n-k)^{-x} \exp (-\lambda_n)\frac{\lambda_n^k}{k!}.
\end{eqnarray*}

\begin{lemma}\label{useddd1}
For all sufficiently large $n$,
$\lambda_n^{-{1}/{2}} ( z_n(x) - \lambda_n  )$ is at least
\[
\exp\bigl(-4 T^{-1}\bigr) \frac{ R_{1,n} }{ R_{2,n} + 4 (1-x)^{-1} T^{-(1-x)} },
\]
and at most
\[
\exp\bigl(4 T^{-1}\bigr) \frac{ R_{1,n} + 4 (1-x)^{-1} T^{-(1-x)} }{ R_{2,n} }.
\]
\end{lemma}
\begin{pf}
The proof is deferred to the \hyperref[2appsec2]{Appendix}.
\end{pf}
Letting $z_{\infty}(x) \stackrel{\Delta}{=} \upsilon(x,-B) + B$, we now
use Lemma~\ref{useddd1} to demonstrate the following:
%
\begin{proposition}\label{astep1}
For $x \in(0,1)$, $\lim_{n \rightarrow\infty} \lambda_n^{-{1}/{2}}  ( z_n(x) - \lambda_n  ) = z_{\infty}(x)$.
\end{proposition}
\begin{pf}
We proceed by relating $R_{1,n}$ and $R_{2,n}$ to the expectations of
certain functions of a scaled Poisson r.v., and then analyze these
expectations as $n \rightarrow\infty$ using tools from
weak-convergence theory. Let $X_n$ denote a Poisson r.v. with mean
$\lambda_n$,
$Z_n \stackrel{\Delta}{=} \lambda_n^{-{1}/{2}} (X_n - \lambda_n)$,
\[
Y_{1,n} \stackrel{\Delta} {=} \biggl( B \biggl(\frac{n}{\lambda_n}
\biggr)^{{1}/{2}} - Z_n \biggr)^{1-x} I \biggl(
Z_n \leq B \biggl(\frac{n}{\lambda_n}\biggr)^{{1}/{2}} - (T+1)
\lambda_n^{-{1}/{2}} \biggr)
\]
and
\begin{eqnarray*}
Y_{2,n}& \stackrel{\Delta} {=}&\biggl(B \biggl(\frac{n}{\lambda_n}
\biggr)^{{1}/{2}} - Z_n \biggr)^{-x} \\
&&{}\times I \biggl(
Z_n \leq\bigl(B - T^{-1}\bigr) \biggl(\frac{n}{\lambda_n}
\biggr)^{{1}/{2}} \\
&&\qquad{}+ \lambda_n^{-{1}/{2}} \bigl( \bigl\lceil n
- T^{-1} n^{{1}/{2}} \bigr\rceil- \bigl(n - T^{-1}
n^{{1}/{2}}\bigr) \bigr) \biggr).
\end{eqnarray*}
It follows from a straightforward computation that $R_{1,n} = \E
[Y_{1,n}]$, and $R_{2,n} = \lambda_n^{-{1}/{2}} \E[ Z_n Y_{2,n}
] +
E[Y_{2,n}]$.
Let $f_1(y) \stackrel{\Delta}{=} (B - y)^{1-x}I(y \leq B), f_2(y)
\stackrel{\Delta}{=}\break (B - y)^{-x}I(y \leq B - T^{-1})$, $f_3(y)
\stackrel{\Delta}{=} y(B - y)^{-x}I(y \leq B - T^{-1})$ and $N$ denote
a normal r.v. with zero mean and unit variance. It may be easily
verified that $\lbrace Y_{1,n} \rbrace, \lbrace Y_{2,n} \rbrace$ and
$\lbrace Z_n Y_{2,n} \rbrace$ are uniformly integrable sequences of
r.v.s, and converge in distribution to $f_1(N), f_2(N), f_3(N)$, respectively. It follows that
$\lim_{n \rightarrow\infty}
\E[ Y_{1,n}]
= \E[ f_1(N) ]
= (2\pi)^{-{1}/{2}}\int_{-\infty}^B (B-y)^{1-x} \exp(-\frac
{y^2}{2}) \,dy$,
$\lim_{n \rightarrow\infty}
\E[ Y_{2,n} ]
=
\E[ f_2(N) ]
=
(2\pi)^{-{1}/{2}}
\int_{-\infty}^{ B - T^{-1} } (B-y)^{-x} \exp( -\frac{y^2}{2} ) \,dy,$
and\break
$\lim_{n \rightarrow\infty}
\E[Z_n Y_{2,n}]
=
\E[f_3(N)]
=
(2\pi)^{-{1}/{2}}
\int_{-\infty}^{ B - T^{-1} } y (B-y)^{-x} \exp( -\frac{y^2}{2} )
\,dy$.\break
Plugging the above limits into Lem\-ma~\ref{useddd1}, and letting $T
\rightarrow\infty$, we conclude that
%
\begin{equation}
\label{limitalmost} \lim_{n \rightarrow\infty} \lambda_n^{-{1}/{2}}
\bigl( z_n(x) - \lambda_n \bigr) = \frac{ \int_{-\infty}^B (B-y)^{1-x} \exp( -{y^2}/{2} ) \,dy } {
\int_{-\infty}^{ B } (B-y)^{-x} \exp( -{y^2}/{2} ) \,dy }.
\end{equation}
We now complete the proof by relating the integrals appearing in (\ref
{limitalmost}) to the parabolic cylinder functions. It is stated in
\cite{G80} that for all $x,z \in\reals$,
%
\begin{equation}
\label{paraprop} D_{x+1}(z) - z D_x(z) + x
D_{x-1}(z) = 0.
\end{equation}
Combining (\ref{gonefishin2}) and (\ref{paraprop}), we find that the
right-hand side of (\ref{limitalmost}) equals
\begin{eqnarray*}
\frac{ \int_{0}^{\infty} y^{1-x} \exp( -{(B-y)^2}/{2} ) \,dy } {
\int_{0}^{ \infty} y^{-x} \exp( -{(B-y)^2}/{2} ) \,dy } &=& \frac
{\Gamma(2-x) {(  ( D_x(-B) + B D_{x-1}(-B)  ) }/{(1-x)}) }{
\Gamma
(1-x) D_{x-1}(-B) } \\
&=& z_{\infty}(x),
\end{eqnarray*}
where the final equality follows from the fact that $\frac{\Gamma
(2-x)}{\Gamma(1-x)} = 1-x$.
\end{pf}
We now complete the proofs of Theorem~\ref{psiasymptotics} and
Corollary~\ref{asscorrect}.
\begin{pf*}{Proof of Theorem~\ref{psiasymptotics} and Corollary \ref
{asscorrect}} Since $\Psi_n(x) = z_n(x) - a_n(x)$, Theorem \ref
{psiasymptotics}
follows from Lemma~\ref{aasymptotics} and Proposition~\ref{astep1}.

We now prove Corollary~\ref{asscorrect}. It follows from the
monotonicity of $z_n$ guaranteed by Lemma~\ref{kijima2}(ii) that
for any sufficiently small positive $\varepsilon$ and all sufficiently
large~$n$, one has
\begin{eqnarray*}
\lambda_n^{-{1}/{2}} \biggl( z_n\biggl(
\frac{B^2}{4} + \varepsilon\biggr) - \lambda_n \biggr) &\leq&
\lambda_n^{-{1}/{2}} \bigl( z_n \bigl(
\bigl(n^{{1}/{2}}-\lambda_n^{{1}/{2}}\bigr)^2
\bigr) - \lambda_n \bigr) \\
&\leq& \lambda_n^{-{1}/{2}}
\biggl( z_n\biggl( \frac{B^2}{4} - \varepsilon\biggr) -
\lambda_n \biggr).
\end{eqnarray*}
Thus by Proposition~\ref{astep1}, for all sufficiently small
$\varepsilon
> 0$,
%
\begin{eqnarray}
\label{zeps1} z_{\infty}\biggl(\frac{B^2}{4}+\varepsilon\biggr) &\leq&
\liminf_{n \rightarrow
\infty} \lambda_n^{-{1}/{2}} \bigl( z_n
\bigl( \bigl(n^{{1}/{2}}-\lambda_n^{{1}/{2}}
\bigr)^2 \bigr) - \lambda_n \bigr)
\nonumber
\\[-8pt]
\\[-8pt]
\nonumber
&\leq&
\limsup_{n
\rightarrow
\infty} \lambda_n^{-{1}/{2}} \bigl( z_n
\bigl( \bigl(n^{{1}/{2}}-\lambda_n^{{1}/{2}}
\bigr)^2 \bigr) - \lambda_n \bigr) \leq z_{\infty}
\biggl(\frac
{B^2}{4} - \varepsilon\biggr).
\end{eqnarray}
We now prove that $z_{\infty}$ is continuous in a neighborhood of
$\frac
{B^2}{4}$, from which we conclude that $\lim_{n \rightarrow\infty}
\lambda_n^{-{1}/{2}}  ( z_n ( (n^{{1}/{2}}-\lambda_n^{{1}/{2}})^2  ) - \lambda_n  ) = z_{\infty}( \frac
{B^2}{4} ).$
Indeed, since $D_x(z) > 0$ for all $z \in\reals$ and $x \leq0$, it
follows that $D_{x-1}(-B) > 0$ for $x \leq1$. The continuity of
$z_{\infty}$ on $(-\infty, 1]$ then follows from the fact that
$D_x(-B)$ is an entire function of $x$~\cite{C04}.

Since $a_n  ( ( n^{{1}/{2}} - \lambda_n^{{1}/{2}}
)^2  ) = (\lambda_n n)^{{1}/{2}}$, we also have that
\[
\lim_{n
\rightarrow\infty}\lambda_n^{-{1}/{2}}  \bigl( a_n \bigl(
\bigl(n^{{1}/{2}}-\lambda_n^{{1}/{2}}\bigr)^2  \bigr) - \lambda_n  \bigr) = \frac
{B}{2}.
\]
Combining the above completes the proof, since $\zeta(B) =
z_{\infty}(\frac{B^2}{4}) - \frac{B}{2}$.
\end{pf*}

\section{\texorpdfstring{Asymptotic analysis of $Z^+(\Psi_n)$}
{Asymptotic analysis of Z+(Psi n)}}\label{1zplussec}
In this section we derive the asymptotics of $Z^+(\Psi_n)$ in the HW
regime. In particular, we prove the following:

\begin{theorem}\label{zzeros}If $B<2$ and $\zeta(B) \leq0$, or $B
\geq
2$, then $\lim_{n \rightarrow\infty} Z^+( \Psi_n ) = Z^+( \Psi_{\infty})$.
\end{theorem}
We first prove some additional properties of $Z^+(\Psi_{\infty})$, namely,
%
\begin{lemma}\label{psiprops22}
If $B < 2$ and $\zeta(B) < 0$, or $B \geq2$, then: $\Psi_{\infty}$ has a unique zero $Z^+(\Psi_{\infty}) \in (0, \min(1, \frac
{B^2}{4})  )$;
$\Psi_{\infty}$ is strictly positive on $[0, Z^+(\Psi_{\infty})
)$;
and $\Psi_{\infty}$ is strictly negative on $ ( Z^+(\Psi_{\infty
}), \min(1, \frac{B^2}{4}) ]$.
Alternatively, if $B < 2$ and $\zeta(B) = 0$, then: $\Psi_{\infty}$ is
strictly positive on $[0, \min( 1, \frac{B^2}{4}) )$, and
$Z^+(\Psi_{\infty}) = \frac{B^2}{4}$.
\end{lemma}
\begin{pf}
We begin by proving that $\Psi_{\infty}$ is continuous and strictly
decreasing on $ [0 , \min(1, \frac{B^2}{4})  ]$. Since $\Psi_{\infty} = z_{\infty} - a_{\infty}$,
it suffices to demonstrate the
continuity and monotonicity of $z_{\infty}$ and $a_{\infty}$
separately. We have already shown that $z_{\infty}$ is continuous on
$(-\infty, 1]$, and it follows from Lemma~\ref{kijima2}(ii) and
Proposition~\ref{astep1} that $z_{\infty}$ is nonincreasing on
$[0,1]$. A straightforward calculation demonstrates that $a_{\infty}$
is continuous and strictly increasing on $[0,\frac{B^2}{4}]$. Combining
the above yields the desired result.

We now treat the case $B < 2$ and $\zeta(B) < 0$, or $B \geq
2$. Note that $\Psi_{\infty}(0) > 0$, since $\Psi_{\infty}(0) =
\upsilon
(0,-B) + B$, and by (\ref{parainteger}), $\upsilon(0,-B) > 0$. Also,
$\Psi_{\infty} (\min(1,\break\frac{B^2}{4}) ) < 0$, which we now
demonstrate by a case analysis. If $B < 2$ and $\zeta(B) < 0$, then
$\min(1, \frac{B^2}{4}) = \frac{B^2}{4}$, and $\Psi_{\infty}(\frac
{B^2}{4}) = \zeta(B) < 0$.
Alternatively, if $B \geq2$, then $\min(1, \frac{B^2}{4}) = 1$. But
$\Psi_{\infty}(1) < 0$, since by (\ref{parainteger}), $\Psi_{\infty
}(1) = -B + \frac{1}{2} (B + (B^2 - 4)^{{1}/{2}} ) < 0.$
Combining the above facts completes the proof. The case $B < 2$ and
$\zeta(B) = 0$ follows similarly, and we omit the details.
\end{pf}
We now complete the proof of Theorem~\ref{zzeros}.
\begin{pf*}{Proof of Theorem~\ref{zzeros}}
We first treat the case $B<2$ and $\zeta(B) < 0$, or $B \geq2$, and
begin by demonstrating that
$\liminf_{n \rightarrow\infty} Z^+( \Psi_n ) \geq Z^+( \Psi_{\infty} )$.
Suppose for contradiction that
$\liminf_{n \rightarrow\infty} Z^+( \Psi_n ) < Z^+( \Psi_{\infty}
)$.
Then it follows from Lemma~\ref{psiprops22} that there exists
$\varepsilon> 0$ s.t. $0 < \liminf_{n\rightarrow\infty} Z^+( \Psi_n ) +
\varepsilon< \min(1, \frac{B^2}{4})$, and $\Psi_{\infty} (
\liminf_{n\rightarrow\infty} Z^+( \Psi_n ) + \varepsilon ) > 0$. Thus by
Theorem~\ref{psiasymptotics}, for all sufficiently large $n$,
$\Psi_{n} ( \liminf_{n\rightarrow\infty} Z^+( \Psi_n ) +
\varepsilon
) > 0$, and by the monotonicity of $\Psi_n$ [see Lemma \ref
{kijima2}(iii)], $\Psi_{n}$ is strictly positive on $
(-\infty,\break \liminf_{n\rightarrow\infty} Z^+( \Psi_n ) + \varepsilon)$.
But by the definition of $\liminf$, this implies the existence of an
infinite strictly increasing sequence of integers $\lbrace n_i \rbrace$
s.t.
$\Psi_{n_i} ( Z^+(\Psi_{n_i} )  ) > 0$ for all $i$. This is a
contradiction, since $\Psi_{n_i} ( Z^+(\Psi_{n_i} )  ) = 0$ for
all $i$, and we conclude that
$\liminf_{n \rightarrow\infty} Z^+( \Psi_n ) \geq Z^+( \Psi_{\infty}
)$. The
proof that $\limsup_{n \rightarrow\infty} Z^+( \Psi_n ) \leq Z^+(
\Psi_{\infty} )$, as well as the proofs for the case
$B<2$ and $\zeta(B) = 0$, follow similarly, and we omit the details.
\end{pf*}

\section{\texorpdfstring{The sign of $\zeta$}{The sign of zeta}}\label{1zetazerosec}
In this section we characterize the sign of $\zeta$ on $(0,2)$, proving
the following:
%
\begin{theorem}\label{bstarbigroot2}
$B^* \in(0,2)$. $\zeta$ is strictly positive on $[0, B^*)$ and
strictly negative on $(B^*, 2]$.
\end{theorem}
 We also complete the proof of Proposition~\ref{Bstarprops}.
Although Theorem~\ref{bstarbigroot2} seems clear from Figure \ref
{figzetab}, the formal proof of this fact is somewhat involved, since
a priori it could be the case that $\zeta$ never actually becomes
strictly negative at $B^*$, or that~$\zeta$ has additional zeros on
$(B^*,2]$. We begin by proving a technical lemma about $\upsilon(x,-B)$.
%
\begin{lemma}\label{parabolicconcave}
For any fixed $B > 0$, $\upsilon(x,-B)$ is a concave function of $x$
on $(0,1)$.
\end{lemma}
\begin{pf}
The proof is deferred to the \hyperref[2appsec2]{Appendix}.
\end{pf}
We now prove some bounds for $\varphi'(B) \stackrel{\Delta}{=} \frac
{d}{dB}\varphi(B)$, when it exists.
%
\begin{lemma}\label{bounddB1}
$\varphi$ is a differentiable function on $(0,2)$, and
\[
\varphi'(B) < \bigl(2 B^{-1} - B\bigr) \varphi(B) -
\frac{B^2}{4} - \varphi^2(B) \leq B^{-2} - 1.
\]
\end{lemma}
\begin{pf}
Note that $\upsilon(x,y)$ is a smooth function of $y$ on $(-\infty
,\infty)$ for any fixed $x \leq1$, and
a smooth function of $x$ on $(-\infty,1]$ for any fixed $y \in\reals$.
Indeed, this follows from
the strict positivity of $D_{x-1}(y)$ for each fixed $x \leq1$, and
the fact that $D_x(y)$ is an entire function of $y$ for each fixed $x$
\cite{E53}, and an entire function of $x$ for each fixed $y$ \cite
{C04}. Thus we may apply the multivariate chain rule to $\varphi$. In
light of (\ref{paraprop}), and the fact (stated in~\cite{G80}) that
for all $x,z \in\reals$,
%
\begin{equation}
\label{paraprop2} \frac{d}{dz} D_{x}(z) + \frac{1}{2} z
D_x(z) - x D_{x-1}(z) = 0,
\end{equation}
it then follows from a straightforward computation that $\varphi$ is
differentiable on $(-\infty, 2]$, and
%
\begin{equation}
\label{mvchain3} \varphi'(B) = \frac{B}{2}\frac{d\upsilon}{dx}
\biggl(\frac{B^2}{4},-B\biggr) - \frac
{B^2}{4} - \varphi^2(B)
- B \varphi(B).
\end{equation}
We now bound $\frac{d\upsilon}{dx}(\frac{B^2}{4},-B)$. The mean value
theorem guarantees the existence of $c \in(0, \frac{B^2}{4})$ s.t.
$\frac{d\upsilon}{dx}(c,-B) = (\frac{B^2}{4})^{-1}  ( \upsilon
(\frac
{B^2}{4},-B) - \upsilon(0,-B)  )$.
In light of Lemma~\ref{parabolicconcave}, we conclude that
%
\begin{equation}
\label{bounddx1}\quad  \frac{d\upsilon}{dx}\biggl(\frac{B^2}{4},-B\biggr) \leq\biggl(
\frac{B^2}{4}\biggr)^{-1} \biggl( \upsilon\biggl(
\frac{B^2}{4},-B\biggr) - \upsilon(0,-B) \biggr) < \frac{4}{B^2}
\varphi(B),
\end{equation}
where the final inequality follows from the fact that $\upsilon(0,-B) >
0$ by (\ref{parainteger}).
Combining (\ref{mvchain3}) and (\ref{bounddx1}) proves the first part
of the lemma. It follows that there exists $x_B \in\reals$ s.t.
$\varphi'(B) \leq(2 B^{-1} - B) x_B - x_B^2 - \frac{B^2}{4}$, which is
at most $B^{-2} - 1$ by elementary calculus. Combining the above
completes the proof.
\end{pf}
We now complete the proof of Theorem~\ref{bstarbigroot2} and
Proposition~\ref{Bstarprops}.
\begin{pf*}{Proof of Theorem~\ref{bstarbigroot2} and Proposition
\ref{Bstarprops}}
We first demonstrate that $\zeta$ is strictly positive on $[0, B^*)$.
Indeed, this follows from (\ref{parainteger}), which implies that
$\zeta
(0) =  ( 2^{{1}/{2}} \int_0^{\infty} \exp( - y^2 ) \,dy  )^{-1}
> 0$.

To complete the proof of Theorem~\ref{bstarbigroot2}, we
will first show that $B^* \in(2^{{1}/{2}}, 2)$, and then apply
Lemma~\ref{bounddB1} to prove that $\zeta'(B) < 0$ on $(B^*,2)$.
We show that $B^* \in(2^{{1}/{2}}, 2)$ in two stages, first
proving that $B^* \in(0,2)$. (\ref{parainteger}) implies that $\zeta
(2) = \frac{-2 e^{-1} }{ e^{-1} } + 1 < 0$. That $B^* \in(0,2)$ then
follows from the fact that $\zeta(0) > 0$, and the intermediate value theorem.

We now demonstrate that $B^* > 2^{{1}/{2}}$. It is proven
in~\cite{D85} Theorem 4.1(i) that
%
\begin{equation}
\label{d85res} \gamma_n \geq\inf_{k \geq1} \bigl(
\lambda_n + \min(k,n) - \lambda_n^{{1}/{2}} \bigl(
\min(k-1, n)^{{1}/{2}} + \min(k, n)^{{1}/{2}} \bigr) \bigr).
\end{equation}
Note that for $1 \leq k \leq n$, $\lambda_n + \min(k,n) - \lambda_n^{{1}/{2}}  ( \min(k-1, n)^{{1}/{2}} + \min(k,
n)^{{1}/{2}}  )$ equals
\[
\bigl( \lambda^{{1}/{2}}_n - k^{{1}/{2}}
\bigr)^2 + \frac{ \lambda_n^{{1}/{2}} }{ k^{{1}/{2}} + (k-1)^{{1}/{2}} } \geq \frac
{1}{2} \biggl(
\frac{\lambda_n}{n}\biggr)^{{1}/{2}}.
\]
For all $k \geq n+1$, the right-hand side of (\ref{d85res}) equals $(n^{{1}/{2}} - \lambda_n^{{1}/{2}})^2.$ Combining the above,
we find that
\[
\gamma_n \geq\min \biggl( \frac{1}{2} \biggl(
\frac{\lambda_n}{n}\biggr)^{{1}/{2}} , \bigl(n^{{1}/{2}} -
\lambda_n^{{1}/{2}}\bigr)^2 \biggr).
\]
Recalling that $\lim_{n \rightarrow\infty} (n^{{1}/{2}} -
\lambda_n^{{1}/{2}})^2 = \frac{B^2}{4}$, it follows that
for any fixed $B < 2^{{1}/{2}}$ and all sufficiently large $n$,
$\gamma_n \geq(n^{{1}/{2}} - \lambda_n^{{1}/{2}})^2$.

Now, suppose for contradiction that $B^* < 2^{{1}/{2}}$.
Then combining Lemma~\ref{bounddB1} with the fact that
by construction $\varphi(B^*) = - \frac{B^*}{2}$, we find that
\[
\label{deriv8} \varphi'\bigl(B^*\bigr) < \biggl(\frac{2}{B^*}
- B^* \biggr) \biggl( - \frac{B^*}{2} \biggr) - \frac
{B^{*2}}{4} - \biggl(
- \frac{B^*}{2} \biggr)^2 = -1.
\]
It follows that $\zeta'(B^*) < 0$, since $\zeta'(B^*) = \varphi'(B^*) +
\frac{1}{2}$, and
there exists $B' \in(0, 2^{{1}/{2}})$ s.t. $\zeta(B') < 0$.
Thus if we define all relevant functions (e.g., $\lambda_n, \Psi_n$) in
terms of $B'$,
Corollary~\ref{asscorrect} implies that $\Psi_n ( (n^{{1}/{2}}
- \lambda_n^{{1}/{2}})^2  ) < 0$ for all sufficiently large\vadjust{\goodbreak}
$n$, and
$\gamma_n < (n^{{1}/{2}} - \lambda_n^{{1}/{2}})^2$ by
Proposition~\ref{crossatleastonce}(i).
But this is a contradiction since we have already shown that $B' <
2^{{1}/{2}}$ implies that $\gamma_n \geq(n^{{1}/{2}} -
\lambda_n^{{1}/{2}})^2$ for all sufficiently large $n$, showing that $B^*
> 2^{{1}/{2}}$.

We now complete the proof of Theorem~\ref{bstarbigroot2} by
demonstrating that $\zeta'(B) < 0$ on $(2^{{1}/{2}},2)$.
Indeed, for $B \in(2^{{1}/{2}},2)$, we have by Lemma \ref
{bounddB1} that $\zeta'(B)$ equals
\[
\varphi'(B) + \frac{1}{2} < \frac{1}{B^2} - 1 +
\frac{1}{2} = 0,
\]
completing the proof of Theorem~\ref{bstarbigroot2}.

We now prove Proposition~\ref{Bstarprops}. In light of
Theorem~\ref{bstarbigroot2}, the value of $B^*$ may easily be
evaluated numerically to the approximate value 1.85772. The second part
of the proposition follows from Lemma~\ref{psiprops22}.
\end{pf*}

\section{Limiting spectral gap in the HW regime and asymptotic phase
transition}\label{1gammansec}
In this section we complete the proofs of Theorem \ref
{asymptoticphasetransition} and Corollary~\ref{asymptoticphasetransitioncor}.
\begin{pf*}{Proof of Theorem~\ref{asymptoticphasetransition}}
First, suppose $0 < B < B^*$. Then it follows from Theorem \ref
{bstarbigroot2} that $B < 2$ and $\zeta(B) > 0$. Combining with
Corollary~\ref{asscorrect}, we conclude that $\Psi_n  (
(n^{{1}/{2}} - \lambda_n^{{1}/{2}})^2  ) > 0$ for all sufficiently
large $n$, and $\gamma_n = (n^{{1}/{2}} - \lambda_n^{{1}/{2}})^2$ by
Proposition~\ref{crossatleastonce}(ii). Observing that
$\lim_{n \rightarrow\infty} (n^{{1}/{2}} - \lambda_n^{{1}/{2}})^2 =
\frac{B^2}{4}$ completes the proof for this case.

Now, suppose $B = B^*$. By Proposition \ref
{crossatleastonce}, for all sufficiently large $n$, either $\Psi_n (
(n^{{1}/{2}}-\lambda_n^{{1}/{2}})^2  ) < 0$, in which case
$\gamma_n = Z^+(\Psi_n)$, or $\gamma_n = (n^{{1}/{2}} - \lambda_n^{{1}/{2}})^2$. Let $\lbrace n_i, i \geq1 \rbrace$ denote the
subsequence of $\lbrace n \rbrace$ for which
$\Psi_{n_i} ((n_i^{{1}/{2}}-\lambda_{n_i}^{{1}/{2}})^2 )
< 0$. If $\lbrace n_i, i \geq1 \rbrace$ is a finite set, then
trivially $\gamma_n = (n^{{1}/{2}} - \lambda_n^{{1}/{2}})^2$
for all sufficiently large $n$, and observing that $\lim_{n
\rightarrow
\infty} (n^{{1}/{2}} - \lambda_n^{{1}/{2}})^2 = \frac
{B^{*2}}{4}$ completes the proof.
Alternatively, suppose $\lbrace n_i, i \geq1 \rbrace$ is an infinite
set. Then Theorem~\ref{zzeros} implies that
$\lim_{i \rightarrow\infty} Z^+( \Psi_{n_i} ) = \frac{B^{*2}}{4}$.
Combining the above completes the proof for this case, since $\gamma_n$
always belongs to one of two series, both of which converge to $\frac
{B^{*2}}{4}$.

Next, consider the case $B \in(B^*, 2)$. It follows from
Theorem~\ref{bstarbigroot2} that $\zeta(B) < 0$. Combining with
Corollary~\ref{asscorrect}, we conclude that $\Psi_n  (
(n^{{1}/{2}} - \lambda_n^{{1}/{2}})^2  ) < 0$ for all sufficiently
large $n$, and $\gamma_n = Z^+(\Psi_n)$ by Proposition \ref
{crossatleastonce}(i). That $\lim_{n \rightarrow\infty}
\gamma_n = Z^+(\Psi_{\infty})$ then follows from Theorem~\ref{zzeros}.

Finally, suppose $B \geq2$. Then $(n^{{1}/{2}} -
\lambda_n^{{1}/{2}})^2 \geq1$ for all sufficiently large $n$, and
Proposition~\ref{crossatleastonce}(iii) implies that $\gamma_n
= Z^+(\Psi_n)$. The proof then follows from Theorem~\ref{zzeros}.
\end{pf*}

\begin{pf*}{Proof of Corollary~\ref{asymptoticphasetransitioncor}}
$\!\!\!$Suppose for contradiction that $\liminf_{n \rightarrow\infty}
n^{{1}/{2}}\times  (1 - \rho^*_n) < B^*$.
Then there exists $\varepsilon> 0$, and an infinite, strictly increasing
sequence of integers $\lbrace n_i, i \geq1 \rbrace$, s.t.
$\rho^*_{n_i} > 1 - (B^* - \varepsilon)n^{-{1}/{2}}_i $ for all $i$.
Consider the sequence $\lbrace Z_i, i \geq1 \rbrace$ of continuous
time Markov chains, in which $Z_i$ is an $M/M/n_i$ queueing system with\vadjust{\goodbreak}
$\lambda_{n_i} = n_i - (B^* - \varepsilon) n^{{1}/{2}}_i$, $\mu= 1$.
Let us define all relevant functions (e.g., $\Psi_{n_i}, \lambda_{n_i}$) w.r.t. $B^* - \varepsilon$. Then since
$\zeta(B^* - \varepsilon) > 0$ and $B^* - \varepsilon< 2$ by Theorem \ref
{bstarbigroot2}, it follows from
Corollary~\ref{asscorrect} that $\Psi_{n_i}  ( ( n^{{1}/{2}}_i
- \lambda^{{1}/{2}}_{n_i} )^2  ) > 0$ for all sufficiently large
$i$, and
$\gamma_{n_i} = ( n^{{1}/{2}}_i - \lambda^{{1}/{2}}_{n_i} )^2$
by Proposition~\ref{crossatleastonce}(ii). But $\frac{
\lambda_{n_i} }{n_i \mu} = 1 - (B^* - \varepsilon)n^{-{1}/{2}}_i < \rho^*_{n_i}$ for all $i$.
This is a contradiction, since by Theorem~\ref{fixedntransition},
$\frac{ \lambda_{n_i} }{n_i \mu} < \rho^*_{n_i} $ implies that
the spectral gap $\gamma_{n_i}$ of $Z_i$ is strictly less than $(
n^{{1}/{2}}_i - \lambda^{{1}/{2}}_{n_i} )^2$. Thus
$\liminf_{n \rightarrow\infty} n^{{1}/{2}} (1 - \rho^*_n) \geq
B^*$. A similar argument demonstrates that $\limsup_{n \rightarrow
\infty} n^{{1}/{2}} (1 - \rho^*_n) \leq B^*$, and we omit the
details. Combining the above completes the proof.
\end{pf*}

\section{Explicit bounds on the distance to stationarity}\label{explicitsec}

In this section we complete the proof of Theorem~\ref{explicitboundstheorem}.

\subsection{KM representation}\label{subseckm}
In this subsection we formally state the KM representation for the
transient distribution of the $M/M/n$ queue, when the traffic intensity
is at least $\rho^*_n$. Let
%
\begin{equation}
\label{qdef} Q_{n,k}(x) \stackrel{\Delta} {=} %
\cases{ 1, &\quad  $
\mbox{if k = 0,}$\vspace*{2pt}
\cr
\displaystyle 1 - \frac{x}{\lambda_n}, &\quad  $\mbox{if k = 1,}$
\vspace*{2pt}
\cr
\displaystyle\biggl(1 - \frac{x}{ \lambda_n } + \frac{ \min( k-1, n) }{ \lambda_n } \biggr)
Q_{n,k-1}(x) \vspace*{2pt}\cr
{}\qquad\displaystyle- \frac{ \min(k-1,n) }{\lambda_n} Q_{n,k-2}(x), & \quad $\mbox{otherwise}$} %
\end{equation}
and
\[
c_n(x) \stackrel{\Delta} {=} Q^2_{n,n}(x) -
\frac{ \lambda_n + n - x
}{\lambda_n } Q_{n,n}(x) Q_{n,n-1}(x) + \frac{n}{ \lambda_n }
Q^2_{n,n-1}(x).
\]
It is proven in~\cite{D81} that $c_n$ is strictly positive on $ ( (
n^{{1}/{2}} - \lambda_n^{{1}/{2}} )^2, ( n^{{1}/{2}} +
\lambda_n^{{1}/{2}} )^2  )$. We also define
\[
b_n(x) \stackrel{\Delta} {=} %
\cases{\displaystyle \bigl(x -
\bigl(n^{{1}/{2}} - \lambda_n^{{1}/{2}}\bigr)^2
\bigr)^{{1}/{2}} \bigl( \bigl(n^{{1}/{2}} + \lambda_n^{{1}/{2}}
\bigr)^2 - x \bigr)^{{1}/{2}}, \vspace*{2pt}\cr
\hspace*{38pt}\mbox{if } \bigl(n^{{1}/{2}}
- \lambda_n^{{1}/{2}}\bigr)^2 \leq x \leq
\bigl(n^{{1}/{2}} + \lambda_n^{{1}/{2}}\bigr)^2,
\vspace*{2pt}
\cr
\infty,\qquad \mbox{otherwise,}} %
\]
and let $g_n(k) \stackrel{\Delta}{=} \lambda_n^{k-n} n^{\min(n-k,0)}
\prod_{i=k+1}^n i$. Then the following is proven by KM in~\cite{KM58}
(see also~\cite{D81}):
%
\begin{theorem}\label{maindoorn}
If $\frac{ \lambda_n }{n} \geq\rho^*_n$, then for all $i,j,t \geq0$,
\begin{eqnarray*}
P^n_{i,j}(t) - P^n_j(\infty) &=& (2
\pi)^{-1}g_n(j) (\lambda_n n)^{-1}
\\
&&{}\times\int_{ ( n^{{1}/{2}} -
\lambda
_n^{{1}/{2}} )^2 }^{ (
n^{{1}/{2}} + \lambda_n^{{1}/{2}} )^2 } \exp(- x t) Q_{n,i}(x)
Q_{n,j}(x) b_n(x) c_n(x)^{-1} \,dx.
\end{eqnarray*}
\end{theorem}

\subsection{Bounds for $ |Q_{n,n}(x) |$, $
|Q_{n,n-1}(x)
|$ and $ |Q_{n,n}(x) - Q_{n,n-1}(x) |$}\label{subsecbase}
In this subsection we prove bounds for $ |Q_{n,n}(x) |$, $
|Q_{n,n-1}(x) |$ and $ |Q_{n,n}(x) -\break  Q_{n,n-1}(x) |$. Let
$h_n(x) \stackrel{\Delta}{=} 2 n b_n(x)^{-1}$. Then we have the following:
%
\begin{lemma}\label{firstbounds}
For all $x \in ( (n^{{1}/{2}} - \lambda_n^{{1}/{2}} )^2,
(n^{{1}/{2}} + \lambda_n^{{1}/{2}} )^2 )$,
%
\begin{eqnarray}
\bigl|Q_{n,n}(x) \bigr| &\leq& c_n(x)^{{1}/{2}}
h_n(x),\label{anothereee1}
\\
\bigl|Q_{n,n-1}(x)\bigr | &\leq& c_n(x)^{{1}/{2}}
h_n(x),\label{anothereee2}
\\
\bigl|Q_{n,n}(x) - Q_{n,n-1}(x)\bigr | &\leq&\biggl(\frac{x}{n}
\biggr)^{{1}/{2}} c_n(x)^{{1}/{2}} h_n(x).
\label{anothereee3}
\end{eqnarray}
\end{lemma}
\begin{pf}
We first prove (\ref{anothereee1}). If $Q_{n,n}(x) = 0$, then $
|Q_{n,n}(x) | = 0 <\break c_n(x)^{{1}/{2}} h_n(x)$.
Otherwise,
%
\begin{eqnarray}
Q^2_{n,n}(x) c_n(x)^{-1} &=& \biggl(
1 - \frac{ \lambda_n + n - x }{ \lambda_n } \frac{
Q_{n,n-1}(x) }{ Q_{n,n}(x) } + \frac{n}{ \lambda_n } \biggl(
\frac{
Q_{n,n-1}(x) }{Q_{n,n}(x) } \biggr)^2 \biggr)^{-1}
\nonumber
\\[-8pt]
\\[-8pt]
\nonumber
&\leq& \sup_{z \in\reals} \biggl( \biggl( 1 - \frac{ \lambda_n + n - x }{
\lambda_n } z +
\frac{n}{ \lambda_n } z^2 \biggr)^{-1} \biggr) = 4
\lambda_n n b_n(x)^{-2},
\end{eqnarray}
where the final equality follows from elementary calculus. Taking
square roots completes the proof. The proof of (\ref{anothereee2})
follows from a similar argument, and we omit the details. We now prove
(\ref{anothereee3}). It is shown in~\cite{D81} that $Q_{n,n}$ and
$Q_{n,n-1}$ do not have any common zeros.
Thus first suppose $Q_{n,n}(x) = 0$. Then $ ( Q_{n,n}(x) -
Q_{n,n-1}(x)  )^2 c_n(x)^{-1} = \frac{ \lambda_n }{n} < 1$.
Combining with the fact that
$4 \lambda_n x b_n(x)^{-2} = 1 + (\lambda_n + x - n)^2 b_n(x)^{-2}
\geq
1$ completes the proof.
The case $Q_{n,n-1}(x) = 0$ follows from a similar argument, and we
omit the details. Finally, suppose $Q_{n,n}(x) \not= 0$ and
$Q_{n,n-1}(x) \not= 0$. Then
%
\begin{eqnarray}\label{usepart23}
&&\bigl( Q_{n,n}(x) - Q_{n,n-1}(x) \bigr)^2
c_n(x)^{-1}\nonumber\\
 &&\qquad= { \bigl(\bigl({ Q_{n,n}(x) }/{Q_{n,n-1}(x)}\bigr) - 1\bigr)^2 }\nonumber\\
 &&\qquad\quad{}/\bigl( \bigl({ Q_{n,n}(x)
}/{Q_{n,n-1}(x)}\bigr)^2 - \bigl({ (\lambda_n + n - x )}/{
\lambda_n }\bigr)\bigl( { Q_{n,n}(x) }/{Q_{n,n-1}(x)}\bigr) \\
&&\hspace*{277pt}\qquad{}+ {n}/{ \lambda_n }\bigr)
\nonumber\\
&&\qquad\leq \sup_{z \in\reals} \biggl( (z-1)^2 \biggl(z^2 -
\frac{ \lambda_n + n - x }{
\lambda_n } z + \frac{n}{ \lambda_n } \biggr)^{-1} \biggr).\nonumber
\end{eqnarray}

Let $f(z) \stackrel{\Delta}{=} (z-1)^2 (z^2 - \frac{ \lambda_n + n
- x
}{ \lambda_n } z + \frac{n}{ \lambda_n } )^{-1}$. It may be easily
verified that $f(z)$ is a continuously differentiable rational function
of $z$ on $\reals$, and the zeros of $\frac{d}{dz} f(z)$ occur at $z =
1$ and $z = \frac{ \lambda_n - n - x }{ \lambda_n - n + x}$. Thus
$\sup_{z \in\reals} f(z)$ must be one of $f(1), f( \frac{ \lambda_n - n
- x
}{ \lambda_n - n +
x} ), \lim_{z \rightarrow- \infty} f(z), \lim_{z \rightarrow\infty
} f(z)$.\vadjust{\goodbreak}
It follows from a straightforward computation that $f(1) = 0$, $\lim_{z
\rightarrow- \infty} f(z) = \lim_{z \rightarrow\infty} f(z) = 1$, and
$f( \frac{ \lambda_n - n - x }{ \lambda_n - n + x } ) = 4 \lambda_n x
b_n(x)^{-2}$.
Combining with (\ref{usepart23}), and the fact that\break $4 \lambda_n x
b_n(x)^{-2} \geq1$, completes the proof.
\end{pf}

\subsection{\texorpdfstring{Bounding $|Q_{n,k}(x) |$ and $ |Q_{n,k \pm 1}(x)
- Q_{n,k}(x) |$ for $k = n \pm O(n^{{1}/{2}})$}
{Bounding $|Q_{n,k}(x) |$ and $|Q_{n,k +- 1}(x)
- Q_{n,k}(x) |$ for $k = n +- O(n^{{1}/{2}})$}}\label{subsecinduct}
In this subsection, we bound $ |Q_{n,k}(x) |$ and $ |Q_{n,k
\pm1}(x) - Q_{n,k}(x) |$ for $k = n \pm O(n^{{1}/{2}})$. Let
$s_n(a) \stackrel{\Delta}{=} a (n^{{1}/{2}} - a)^{-1}$, $r_n(a,x)
\stackrel{\Delta}{=} x (n - a n^{{1}/{2}})^{-1}$ and
$F_n(a,x) \stackrel{\Delta}{=} \exp (  ( 1 + s_n(a)  )
(a
+ n^{-{1}/{2}} ) ( 3 x^{{1}/{2}} + a )  )$.
Then we
prove the following:
%
\begin{theorem}\label{fundlemmabound1}
For all $a \geq B > 0$, $k \in[n - a n^{{1}/{2}}, n + a n^{{1}/{2}} + 1]$ and
$x \in ( (n^{{1}/{2}} - \lambda_n^{{1}/{2}} )^2,
(n^{{1}/{2}} + \lambda_n^{{1}/{2}} )^2 )$,
%
\begin{eqnarray}
\bigl|Q_{n,k}(x)\bigr | &\leq& F_n(a,x) c_n(x)^{{1}/{2}}
h_n(x),
\\
\bigl|Q_{n,k+1}(x) - Q_{n,k}(x)\bigr | &\leq & r_n(a,x)^{{1}/{2}}
F_n(a,x) c_n(x)^{{1}/{2}} h_n(x).
\end{eqnarray}
\end{theorem}

We first bound $ |Q_{n, k \pm1}(x) |$ and $ |Q_{n,k \pm1}(x)
- Q_{n,k}(x) |$ in terms of $ |Q_{n,k}(x) |$ and $
|Q_{n,k}(x) - Q_{n,k \pm1}(x) |$. Namely,
%
\begin{lemma}\label{moreqbounds1}
For all $a \geq B > 0$, $k \geq n - a n^{{1}/{2}}$, $x > 0$ and $i
\in\lbrace1,-1\rbrace$,
%
\begin{eqnarray}\label{provefirst1}
&&\bigl|Q_{n,k+i}(x) \bigr|
\nonumber
\\[-8pt]
\\[-8pt]
\nonumber
&&\qquad\leq\exp \bigl( r_n(a,x) +
s_n(a) \bigr) \bigl( \bigl|Q_{n,k}(x) \bigr| + \bigl|Q_{n,k}(x) -
Q_{n,k-i}(x)\bigr | \bigr),
\\
&&\bigl|Q_{n,k+i}(x) - Q_{n,k}(x)\bigr |
\nonumber
\\[-8pt]
\\[-8pt]
\nonumber
\qquad&&\qquad\leq\exp \bigl(
r_n(a,x) + s_n(a) \bigr) \bigl( r_n(a,x)
\bigl|Q_{n,k}(x)\bigr |
+\bigl |Q_{n,k}(x) - Q_{n,k-i}(x)\bigr | \bigr).
\label{provefirst2}
\end{eqnarray}
\end{lemma}

\begin{pf}
Note that
\begin{eqnarray*}
\bigl|Q_{n,k+i}(x) \bigr| &= &\biggl| \bigl(1 - x \min(k,n)^{{(i-1)}/{2}}
\lambda_n^{{(-i-1)}/{2}} \bigr) Q_{n,k}(x)\\
&&\hspace*{8pt}{} + \biggl(
\frac{ \min( k,
n) }{
\lambda_n } \biggr)^i \bigl( Q_{n,k}(x) -
Q_{n,k-i}(x) \bigr) \biggr|,
\\
\bigl|Q_{n,k+i}(x) - Q_{n,k}(x) \bigr| &=& \biggl|- \frac{x}{ \lambda_n }
Q_{n,k}(x) + \biggl(\frac{ \min( k, n) }{ \lambda_n } \biggr)^i \bigl(
Q_{n,k}(x) - Q_{n,k-i}(x) \bigr) \biggr|.
\end{eqnarray*}
Since $\max (  |1 - \frac{x}{\lambda_n} |,  |1 -
\frac
{x}{\min(k,n)} | ) \leq\exp (r_n(a,x) )$,
$\max (\frac{\min(k,n)}{\lambda_n},\frac{\lambda_n}{\min
(k,n)} )
\leq\exp (s_n(a) )$
and $|\frac{x}{\lambda_n}| \leq r_n(a,x)$, the proof then follows from
the triangle inequality.
\end{pf}
 We now use an induction argument to bound $
|Q_{n,k}(x)
|$ and $ |Q_{n,k}(x) - Q_{n,k \pm1}(x) |$ for $k = n \pm
O(n^{{1}/{2}})$. Let $G_n(a,x) \stackrel{\Delta}{=} \exp
(r_n(a,x) + r_n(a,x)^{{1}/{2}} + s_n(a) )$. Then we demonstrate
the following:
%
\begin{lemma}\label{boundgrowthqprelim}
For all $a \geq B > 0$, $k \geq n - a n^{{1}/{2}}$ and $x > 0$,
%
\begin{eqnarray}
\bigl|Q_{n,k}(x)\bigr | &\leq &G_n(a,x)^{|k-n|}
c_n(x)^{{1}/{2}} h_n(x),\label{assertme1}
\\
\qquad\bigl|Q_{n,k}(x) - Q_{n, k + 1 - 2 I(k \geq n)}(x) \bigr| &\leq& r_n(a,x)^{{1}/{2}}
G_n(a,x)^{|k-n|} c_n(x)^{{1}/{2}}
h_n(x).\label{assertme2}
\end{eqnarray}
\end{lemma}
\begin{pf}
We first treat the case $k \geq n$. We proceed by induction on (\ref
{assertme1}) and~(\ref{assertme2}) simultaneously. The base case $k =
n$ follows immediately from Lemma~\ref{firstbounds}. Now, suppose the
induction is true for some $k \geq n$.
Then by Lemma~\ref{moreqbounds1} and the induction hypothesis $
|Q_{n,k+1}(x) |$ is at most
\begin{eqnarray*}
& & \exp \bigl( r_n(a,x) + s_n(a) \bigr)\\
&&\quad{}\times \bigl(
G_n(a,x)^{k-n} c_n(x)^{{1}/{2}}
h_n(x) + r_n(a,x)^{{1}/{2}} G_n(a,x)^{k-n}
c_n(x)^{{1}/{2}} h_n(x) \bigr)
\\
& & \quad= \exp \bigl( r_n(a,x) + s_n(a) \bigr)
G_n(a,x)^{k-n} c_n(x)^{{1}/{2}}
h_n(x) \bigl(1 + r_n(a,x)^{{1}/{2}} \bigr)
\\
& & \quad \leq G_n(a,x)^{k + 1 - n} c_n(x)^{{1}/{2}}
h_n(x).
\end{eqnarray*}
Similarly, by Lemma~\ref{moreqbounds1} and the induction hypothesis,
$ |Q_{n,k+1}(x) - Q_{n,k}(x) |$ is at most
\begin{eqnarray*}
& &\exp \bigl( r_n(a,x) + s_n(a) \bigr)\\
&&\quad{}\times \bigl(
r_n(a,x) G_n(a,x)^{k-n} c_n(x)^{{1}/{2}}
h_n(x) \\
&&\hspace*{6pt}\qquad{}+ r_n(a,x)^{{1}/{2}} G_n(a,x)^{k-n}
c_n(x)^{{1}/{2}} h_n(x) \bigr)
\\
& & \qquad = \exp \bigl( r_n(a,x) + s_n(a) \bigr)
G_n(a,x)^{k-n} c_n(x)^{{1}/{2}}
h_n(x)\\
&&\hspace*{30pt}{}\times \bigl(1 + r_n(a,x)^{{1}/{2}} \bigr)
r_n(a,x)^{{1}/{2}}
\\
& &\qquad  \leq r_n(a,x)^{{1}/{2}} G_n(a,x)^{k + 1 - n}
c_n(x)^{{1}/{2}} h_n(x).
\end{eqnarray*}
This concludes the induction, proving (\ref{assertme1}) and (\ref
{assertme2}) for the case $k \geq n$.

The proof for the case $k < n$ follows from a similar
argument, and we omit the details.
\end{pf}
With Lemma~\ref{boundgrowthqprelim} in hand, we now complete the proof
of Theorem~\ref{fundlemmabound1}.
\begin{pf*}{Proof of Theorem~\ref{fundlemmabound1}}
By Lemma~\ref{boundgrowthqprelim}, $ |Q_{n,k}(x) |$ is at most
%
\begin{eqnarray}\label{expoin11}
& & \exp \bigl( \bigl(a n^{{1}/{2}} + 1 \bigr) \bigl( r_n(a,x)
+ r_n(a,x)^{{1}/{2}} + s_n(a) \bigr) \bigr)
c_n(x)^{{1}/{2}} h_n(x)
\nonumber
\\
& & \quad = \exp \bigl( \bigl(a + n^{-{1}/{2}} \bigr) \bigl(1 + s_n(a)
\bigr)
\nonumber
\\[-8pt]
\\[-8pt]
\nonumber
&&\hspace*{38pt}{}\times \bigl( x n^{-{1}/{2}} + \bigl(1 + s_n(a)
\bigr)^{-{1}/{2}} x^{{1}/{2}} + a \bigr) \bigr) c_n(x)^{{1}/{2}}
h_n(x)
\\
& &\quad  \leq \exp \bigl( \bigl(a + n^{-{1}/{2}} \bigr) \bigl(1 +
s_n(a) \bigr) \bigl( x n^{-{1}/{2}} + x^{{1}/{2}} + a \bigr)
\bigr) c_n(x)^{{1}/{2}} h_n(x).\nonumber
\end{eqnarray}
Similarly, $ |Q_{n,k+1}(x) - Q_{n,k}(x) |$ is at most
%
\begin{equation}
\label{ineqexpo1b}r_n(a,x)^{{1}/{2}} \exp \bigl( \bigl(a +
n^{-{1}/{2}} \bigr) \bigl(1 + s_n(a) \bigr) \bigl( x
n^{-{1}/{2}} + x^{{1}/{2}} + a \bigr) \bigr) c_n(x)^{{1}/{2}}
h_n(x).\hspace*{-35pt}\vadjust{\goodbreak}
\end{equation}
Furthermore, note that $x n^{-{1}/{2}} < 2 x^{{1}/{2}}$ for $x
\in (0, (n^{{1}/{2}} + \lambda_n^{{1}/{2}})^2 )$, since
\[
\frac{ x n^{-{1}/{2}} }{2x^{{1}/{2}}} = \frac{ x^{{1}/{2}} }{ 2 n^{{1}/{2}} } < \frac{ n^{{1}/{2}} + \lambda_n^{{1}/{2}} }{ 2 n^{{1}/{2}} } < 1.
\]
Combining with (\ref{expoin11}) and (\ref{ineqexpo1b}) completes the proof.
\end{pf*}

\subsection{\texorpdfstring{Proof of Theorem \protect\ref{explicitboundstheorem}}
{Proof of Theorem 5}}\label{subsecproof}
In this subsection we complete the proof of Theorem \ref
{explicitboundstheorem}.
 We begin by deriving a variant of the KM representation for
$P^n_{i,\leq j}(t)$, as opposed to
$P^n_{i,j}(t)$, that does not simply sum over all $j+1$ states $\leq
j$, but instead relies on a ``probability flow'' interpretation using the
Chapman--Kolmogorov (CK) differential equations.
%
\begin{lemma} \label{pleq}
If $\frac{ \lambda_n }{n} \geq\rho^*_n$, then for all $i,j,t \geq0$,
$|P^n_{i,\leq j}(t) - P^n_{\leq j}(\infty)|$ is at most
\begin{eqnarray*}
&&(2 \pi)^{-1} g_n(j) n^{-1} \int
_{ ( n^{{1}/{2}}
- \lambda_n^{{1}/{2}} )^2 }^{ ( n^{{1}/{2}} + \lambda
_n^{{1}/{2}} )^2 } \exp(- x t) x^{-1} \bigl|
Q_{n,i}(x) \bigr| \\
&&\hspace*{132pt}{}\times\bigl| Q_{n,j+1}(x) - Q_{n,j}(x)\bigr |
b_n(x) c_n(x)^{-1} \,dx.
\end{eqnarray*}
\end{lemma}

\begin{pf}
The CK differential equations imply that $\frac{d}{dt} P^n_{i, \leq j
}(t) = \min(j+1,n) P^n_{i,j+1}(t) - \lambda_n P^n_{i,j}(t)$. Thus for
all $i,j,t \geq0$,
%
\begin{equation}
\label{labeq11}\quad  \bigl|P^n_{i, \leq j}(t) - P^n_{ \leq j}(
\infty)\bigr | =\biggl |\int_t^{\infty} \bigl( \min(j+1,n)
P^n_{i,j+1}(s) - \lambda_n P^n_{i,j}(s)
\bigr) \,ds \biggr|.
\end{equation}
By detailed balance,
\[
\min(j+1,n) P^n_{j+1}(\infty) - \lambda_n
P^n_{j}(\infty) = 0.
\]
Similarly,
\[
\min(j+1,n) g_n(j+1) = \lambda_n g_n(j).
\]
It thus follows from Theorem~\ref{maindoorn} that the right-hand side
of (\ref{labeq11}) equals
\begin{eqnarray*}
& & \biggl| \int_t^{\infty} \biggl( (2\pi)^{-1}g_n(j)
n^{-1}\\
&&\hspace*{26pt}{}\times \int_{ ( n^{{1}/{2}} -
\lambda_n^{{1}/{2}} )^2 }^{ ( n^{{1}/{2}} + \lambda
_n^{{1}/{2}} )^2 } \exp(- x s)
Q_{n,i}(x)\\
&&\hspace*{94pt}{}\times \bigl( Q_{n,j+1}(x) - Q_{n,j}(x) \bigr)
b_n(x) c_n(x)^{-1} \,dx \biggr) \,ds\biggr |
\\
& &\qquad \leq \int_t^{\infty} \biggl( (2
\pi)^{-1} g_n(j) n^{-1}\\
&&\hspace*{56pt}{}\times \int_{ ( n^{{1}/{2}} - \lambda
_n^{{1}/{2}} )^2 }^{ ( n^{{1}/{2}} + \lambda_n^{{1}/{2}}
)^2 }
\exp(- x s) \bigl|Q_{n,i}(x)\bigr |\\
&&\hspace*{124pt}{}\times \bigl|Q_{n,j+1}(x) - Q_{n,j}(x)
\bigr| b_n(x) c_n(x)^{-1} \,dx \biggr) \,ds
\\
& &\qquad  = (2 \pi)^{-1} g_n(j) n^{-1} \\
&&\hspace*{32pt}{}\times\int
_{ ( n^{{1}/{2}} - \lambda_n^{{1}/{2}} )^2 }^{ (
n^{{1}/{2}} + \lambda_n^{{1}/{2}} )^2 } \exp(- x t) x^{-1}
\bigl|Q_{n,i}(x) \bigr|\\
&&\hspace*{100pt}{}\times \bigl|Q_{n,j+1}(x) - Q_{n,j}(x) \bigr|
b_n(x) c_n(x)^{-1} \,dx,
\end{eqnarray*}
where the final equality follows from Tonelli's theorem. Combining the
above completes the proof.
\end{pf}
 We now prove bounds on a special type of integral that arises
in the analysis of $P^n_{i,j}(t) - P^n_j(\infty)$.
%
\begin{lemma}\label{interintegralbounds}
For all $B,a > 0$ there exists $N_{B,a},C_{B,a} < \infty$, depending
only on $B$ and $a$, s.t. for all $n \geq N_{B,a}$ and $t \geq1$,
%
\begin{eqnarray}
\label{anotherinteqq} &&\int_{(n^{{1}/{2}} - \lambda_n^{{1}/{2}})^2}^{(n^{{1}/{2}}
+ \lambda_n^{{1}/{2}})^2 } \exp(-xt)
F_n(a,x)^2 b_n(x)^{-1} \,dx
\nonumber
\\[-8pt]
\\[-8pt]
\nonumber
&&\qquad\leq
\bigl(1 + C_{B,a} n^{-{1}/{2}}\bigr) \biggl(\frac{\pi}{t \lambda_n}
\biggr)^{{1}/{2}} \exp \biggl( 20 a^2 + 3 a B - \frac{B^2}{4}
t \biggr).
\end{eqnarray}
\end{lemma}
\begin{pf}
The proof is deferred to the \hyperref[2appsec2]{Appendix}.
\end{pf}
Finally, we complete the proof of Theorem~\ref{explicitboundstheorem}.
\begin{pf*}{Proof of Theorem~\ref{explicitboundstheorem}}
Suppose $B \in(0,B^*)$, and $a_1,a_2 \in\reals$. Let $a = \max
(B,|a_1|,|a_2|)$, $i = \lceil n + a_1 n^{{1}/{2}} \rceil$, and $j =
\lceil n + a_2 n^{{1}/{2}} \rceil$. We first prove (\ref
{slayme11}). It follows from Theorem~\ref{maindoorn} and Corollary~\ref{asymptoticphasetransitioncor} that for all sufficiently large $n$
and all $t \geq1$, the left-hand side of (\ref{slayme11}) is at most
%
\begin{eqnarray}
\label{inter22inter} &&(2\pi)^{-1} g_n(j) \bigl(
\lambda_n n^{{1}/{2}}\bigr)^{-1}
\nonumber
\\[-8pt]
\\[-8pt]
\nonumber
&&\qquad{}\times \int
_{ (
n^{{1}/{2}} - \lambda_n^{{1}/{2}} )^2 }^{ (
n^{{1}/{2}} + \lambda_n^{{1}/{2}} )^2 } \exp (- x t ) \bigl|Q_{n,i}(x) \bigr|
\bigl|Q_{n,j}(x)\bigr | c_n(x)^{-1} b_n(x) \,dx.
\end{eqnarray}
Applying Theorem~\ref{fundlemmabound1} to $|Q_{n,i}(x)|$ and
$|Q_{n,j}(x)|$ in (\ref{inter22inter}), we find that the left-hand side
of (\ref{slayme11}) is at most
%
\begin{equation}
\label{anotherinter1} 2 \pi^{-1} g_n(j) \frac{ n^{{3}/{2}} }{\lambda_n}
\int_{ ( n^{{1}/{2}} - \lambda_n^{{1}/{2}} )^2 }^{ (
n^{{1}/{2}} + \lambda_n^{{1}/{2}} )^2 } \exp(- x t ) F_n(a,x)^2
b_n(x)^{-1} \,dx.
\end{equation}
It then follows from Lemma~\ref{interintegralbounds} that there
exists $N_{B,a},C_{B,a} < \infty$, depending only on $B$ and $a$, s.t.
for all $n \geq N_{B,a}$ and $t \geq1$, the left-hand side of (\ref
{slayme11}) is at most
%
\begin{equation}
\label{anotherinterrrr1}\quad  2 (\pi t)^{-{1}/{2}} g_n(j) \biggl(
\frac{n}{\lambda_n}\biggr)^{{3}/{2}} \bigl(1 + C_{B,a}
n^{-{1}/{2}}\bigr) \exp \biggl( 20 a^2 + 3 a B - \frac{B^2}{4}
t \biggr).
\end{equation}
Since $g_n(j) \leq(\frac{n}{\lambda_n})^{n-\lambda_n+1}$, combining
(\ref{anotherinterrrr1}) with a simple Taylor series expansion, and the
fact that $B < B^* < 2$, completes the proof of (\ref{slayme11}).

We now prove (\ref{slayme22}). It follows from Lemma \ref
{pleq} and Corollary~\ref{asymptoticphasetransitioncor} that
for all sufficiently large $n$ and all $t \geq1$, the left-hand side
of (\ref{slayme22}) is at most
%
\begin{eqnarray}
\label{interme1} &&(2\pi)^{-1} g_n(j) n^{-1}\nonumber\\
&&\qquad{}\times \int
_{ ( n^{{1}/{2}}
- \lambda_n^{{1}/{2}} )^2 }^{ ( n^{{1}/{2}} + \lambda
_n^{{1}/{2}} )^2 } \exp(- x t) x^{-1} \bigl|
Q_{n,i}(x) \bigr|\\
&&\hspace*{68pt}\qquad{}\times \bigl| Q_{n,j+1}(x) - Q_{n,j}(x)\bigr |
c_n(x)^{-1} b_n(x) \,dx.\nonumber
\end{eqnarray}
Applying Theorem~\ref{fundlemmabound1} to $ |Q_{n,i}(x) |$ and
$ |Q_{n,j+1}(x) - Q_{n,j}(x) |$, we find that (\ref{interme1}) is
at most
\[
2 \pi^{-1} g_n(j) n \bigl(n-a n^{{1}/{2}}
\bigr)^{-{1}/{2}} \int_{ ( n^{{1}/{2}} - \lambda_n^{{1}/{2}} )^2 }^{ (
n^{{1}/{2}} + \lambda_n^{{1}/{2}} )^2 } \exp(- x t)
x^{-{1}/{2}} F_n(a,x)^2 b_n(x)^{-1}
\,dx.
\]
Since $x^{-{1}/{2}} \leq2 B^{-1}$ for $x \geq( n^{{1}/{2}} -
\lambda_n^{{1}/{2}} )^2$, the proof of
(\ref{slayme22}) then follows from an argument similar to that used to
prove (\ref{slayme11}), and we omit the details.
\end{pf*}

\section{Comparison to other bounds from the literature}\label{comparesec}
In this subsection we compare our bounds from Theorem \ref
{explicitboundstheorem} to two other explicit bounds given in the
literature~\cite{Z91,C98}.
In both cases we will prove that the bounds from the literature
(applied to $ | P^n_{n,\leq n}(t) - P^n_{\leq n}(\infty)  |$ for
$0 < B < B^*$) scale unfavorably in the HW regime.
We begin with the bounds given in~\cite{Z91}, which prove that for
each $B \in(0, B^*)$, there exists $N_B$ s.t. for all $n \geq N_B$ and
$t \geq0$, $ | P^n_{n,\leq n}(t) - P^n_{\leq n}(\infty)  |$ is
at most
%
\begin{eqnarray}
\label{zeifbound}&& 4 (n-1) \Biggl( \sum_{i=1}^{\infty}
\biggl( \biggl(\frac{n}{n-1}\biggr)^i - 1 \biggr)
P^n_i(\infty) + \biggl( \biggl(\frac{n}{n-1}
\biggr)^n - 1 \biggr) \bigl( 1 - 2 P^n_n(
\infty) \bigr) \Biggr)
\nonumber\hspace*{-15pt}
\\[-8pt]
\\[-8pt]
\nonumber
&&\qquad{}\times \exp \bigl( - \bigl(Bn^{{1}/{2}}-1\bigr) (n -
1)^{-1} t \bigr).\hspace*{-15pt}
\end{eqnarray}
Since $\lim_{n \rightarrow\infty}  ( (Bn^{{1}/{2}}-1)(n -
1)^{-1}  ) = 0$, the exponential rate of convergence demonstrated by
(\ref{zeifbound}) goes to zero as $n \rightarrow\infty$, rendering the
bound in~\cite{Z91} ineffective.
We now examine the bounds given in~\cite{C98}, which prove that
for each $B \in(0, B^*)$, there exists\vadjust{\goodbreak} $N_B$ s.t. for all $n \geq N_B$
and $t \geq0$, $ | P^n_{n,\leq n}(t) - P^n_{\leq n}(\infty)  |$
is at most
%
\begin{equation}
\label{chenbound} \bigl( P^n_n(\infty)^{-1} - 1
\bigr)^{{1}/{2}} \exp( - \gamma_n t ).
\end{equation}
It is well known (see~\cite{HW81}) that $\liminf_{n \rightarrow
\infty
}  ( P^n_n(\infty)^{-1} - 1  )^{{1}/{2}} n^{-{1}/{4}} > 0$.
It follows that the prefactor appearing in (\ref{chenbound}) diverges
as $n \rightarrow\infty$, rendering the bound in~\cite{C98} ineffective.

It should be noted that although the bounds given in \cite
{Z91} and~\cite{C98} are ineffective in the HW regime, both bounds
hold in much greater generality, and thus remain interesting and
applicable in a variety of other settings.

\section{Conclusion and open questions}\label{2conc}
In this paper we proved several results about the rate of convergence
to stationarity, that is, the spectral gap, for the $M/M/n$ queue in
the HW regime. We identified the
limiting rate of convergence to steady-state, and proved that an
asymptotic phase transition occurs w.r.t. this rate. In particular, we
demonstrated the
existence of a constant $B^* \approx1.85772$ s.t. when a certain
excess parameter $B \in(0,B^*)$, the error in the steady-state approximation
converges exponentially fast to zero at rate $\frac{B^2}{4}$. For $B >
B^*$, the error in the steady-state approximation converges
exponentially fast to zero at a
different rate, which is the solution to an explicit equation given in
terms of the parabolic cylinder functions. This result may be
interpreted as an
asymptotic version of a phase transition proven to occur for any fixed
$n$ by van Doorn in~\cite{D81}. We also proved explicit bounds on the
distance to stationarity for the $M/M/n$ queue in the HW regime, when
$B < B^*$. Our bounds scale independently of $n$ in the HW regime, and
do not follow from the weak-convergence theory.

This work leaves several interesting directions for future
research. There are many open questions related to the
interaction between weak convergence and convergence to stationarity.
Although our results and those of~\cite{KL08} show that for the
$M/M/n$ queue in the HW regime there is an ``interchange of limits'' in
this regard, namely the limiting rate of convergence equals the rate of
convergence
of the limit, it is unknown to what extent such an interchange must
hold in general. Similarly, it is an open challenge to derive uniform
bounds on the
distance to steady-state in the HW regime for the case of non-Markovian
processing times. It would also be interesting to prove that a phase
transition occurs in other related models, and we refer the reader to
the recent paper~\cite{KL10aa} for some results in this
direction.

\begin{appendix}\label{2appsec2}
\section*{Appendix}
\begin{pf*}{Proof of Lemma~\ref{useddd1}}
By Lemma~\ref{kijima2}(i), $f_{n,n-1}(x)> 0$, and thus $ (
z_n(x) - \lambda_n  )\lambda_n^{-{1}/{2}}$ equals
%
\begin{eqnarray}\label{exp1}
& & \biggl(\frac{\sum_{k=0}^n {n \choose k} \lambda_n^k \prod_{j=1}^{n-k}
(j-x) }{\sum_{k=0}^{n-1} {n-1 \choose k} \lambda_n^k \prod_{j=1}^{n-1-k} (j-x) } - \lambda_n \biggr)
\lambda_n^{-{1}/{2}}
\nonumber
\\[-8pt]
\\[-8pt]
\nonumber
& &\qquad = \frac{\sum_{k=0}^n {n \choose k} \lambda_n^k \prod_{j=1}^{n-k} (j-x)
- \lambda_n \sum_{k=0}^{n-1} {n-1 \choose k} \lambda_n^k \prod_{j=1}^{n-1-k} (j-x) } {
\lambda_n^{{1}/{2}} \sum_{k=0}^{n-1} {n-1 \choose k} \lambda_n^k
\prod_{j=1}^{n-1-k} (j-x) }.
\end{eqnarray}
Note that the numerator of (\ref{exp1}) equals
%
\begin{eqnarray}\label{numerator1}
& &\prod_{j=1}^{n} (j-x) + \sum
_{k=1}^n \pmatrix{n \cr k} \lambda_n^k
\prod_{j=1}^{n-k} (j-x) \nonumber\\
&&\quad{}- \sum
_{k=0}^{n-1} \pmatrix{n-1 \cr (k+1) - 1}
\lambda_n^{k+1} \prod_{j=1}^{n-(k+1)}
(j-x)
\\
& & \qquad= (n-1)! \sum_{k=0}^n (n-k) \prod
_{j=1}^{n-k} \biggl(1-\frac
{x}{j}
\biggr) \frac{ \lambda_n^k }{k!},\nonumber
\end{eqnarray}
and the denominator of (\ref{exp1}) equals
%
\begin{eqnarray}\label{denominator1}
&&\lambda_n^{-{1}/{2}} \sum_{k=0}^{n-1}
\pmatrix{n-1 \cr(k+1)-1} \lambda_n^{k+1} \prod
_{j=1}^{n-(k+1)} (j-x)
\nonumber\\
&&\qquad= \lambda_n^{-{1}/{2}}
\sum_{k=1}^{n} \frac{k}{n}\pmatrix {n
\cr k} \lambda_n^{k} \prod
_{j=1}^{n-k} (j-x)
\\
&&\qquad= \lambda_n^{-{1}/{2}} (n-1)! \sum_{k=0}^{n}
k \prod_{j=1}^{n-k} \biggl(1-
\frac{x}{j}\biggr) \frac{ \lambda_n^k }{k!}.\nonumber
\end{eqnarray}
Plugging (\ref{numerator1}) and (\ref{denominator1}) into (\ref{exp1}),
and multiplying through by $\frac{\exp(-\lambda_n)}{(n-1)!}$, we find
that $ ( z_n(x) - \lambda_n  )\lambda_n^{-{1}/{2}}$ equals
%
\begin{equation}
\label{lastline1} \frac{
\sum_{k=0}^n (n-k) \prod_{j=1}^{n-k} (1-({x}/{j})) \exp(-\lambda_n)({ \lambda_n^k }/{k!})
} {
\lambda_n^{-{1}/{2}} \sum_{k=0}^{n} k \prod_{j=1}^{n-k}
(1-({x}/{j}))   \exp(-\lambda_n)({\lambda_n^k }/{k!}) }.
\end{equation}
We now demonstrate that for all sufficiently large $n$, the numerator
of (\ref{lastline1}) is
at least
\[
\prod_{j=1}^T\biggl(1-\frac{x}{j}
\biggr) \sum_{k=0}^{n-(T+1)} (n - k) \prod
_{j=T+1}^{n-k}\biggl(1-\frac{x}{j}
\biggr) \exp(-\lambda_n)\frac
{\lambda_n^k}{k!},
\]
and at most
\[
\prod_{j=1}^T\biggl(1-\frac{x}{j}
\biggr) \sum_{k=0}^{n-(T+1)} (n - k) \prod
_{j=T+1}^{n-k}\biggl(1-\frac{x}{j}
\biggr) \exp(-\lambda_n)\frac{\lambda_n^k}{k!} + (T+1)^2.
\]
The numerator of (\ref{lastline1}) equals
%
\begin{eqnarray}
\label{numerator2} &&\prod_{j=1}^T\biggl(1-
\frac{x}{j}\biggr) \sum_{k=0}^{n-(T+1)} (n
- k) \prod_{j=T+1}^{n-k}\biggl(1-
\frac{x}{j}\biggr) \exp(-\lambda_n)\frac
{\lambda_n^k}{k!}
\nonumber
\\[-8pt]
\\[-8pt]
\nonumber
&&\qquad{}+ \sum
_{k=n-T}^n (n-k)\prod
_{j=1}^{n-k}\biggl(1-\frac{x}{j}\biggr) \exp(-
\lambda_n)\frac
{\lambda_n^k}{k!}.
\end{eqnarray}
The desired lower bound follows from the fact that the second summand
in (\ref{numerator2}) is nonnegative. The upper bound follows from the
fact that
$\exp(-\lambda_n)\frac{\lambda_n^k}{k!} \leq n^{-{1}/{2}}$ for all
$k \geq0$ by Stirling's inequality, $n - k \leq T + 1$ for all $k \geq
n-T$, and $1 - \frac{x}{j} \leq1$ for all $j \geq1$.

It follows from a similar argument that for all sufficiently
large $n$, the denominator of (\ref{lastline1}) is at least
\[
\lambda_n^{-{1}/{2}} \prod_{j=1}^T
\biggl(1-\frac{x}{j}\biggr) \sum_{k=0}^{n-(T+1)}
k \prod_{j=T+1}^{n-k}\biggl(1-
\frac{x}{j}\biggr) \exp (-\lambda_n)\frac{\lambda_n^k}{k!},
\]
and at most
\[
\lambda_n^{-{1}/{2}} \prod_{j=1}^T
\biggl(1-\frac{x}{j}\biggr) \sum_{k=0}^{n-(T+1)}
k \prod_{j=T+1}^{n-k}\biggl(1-
\frac{x}{j}\biggr) \exp (-\lambda_n)\frac{\lambda_n^k}{k!} +
(T+1)^2,
\]
and we omit the details. Combining the above upper and lower bounds for
the numerator and denominator of (\ref{lastline1}), and dividing
through by
$\prod_{j=1}^T(1-\frac{x}{j})$, we find that for all sufficiently large
$n$, $ ( z_n(x) - \lambda_n  )\lambda_n^{-{1}/{2}}$ is
at least
%
\begin{eqnarray}
\label{usedb1} &&\Biggl( \sum_{k=0}^{n-(T+1)}
(n - k) \prod_{j=T+1}^{n-k}\bigl(1-({x}/{j})\bigr) \exp(-\lambda_n)
\bigl({\lambda_n^k}/{k!}\bigr)\Biggr) \nonumber\\
&&\qquad\Big/\Biggl (
\lambda_n^{-{1}/{2}} \sum_{k=0}^{n-(T+1)} k \prod_{j=T+1}^{n-k}\bigl(1-({x}/{j})\bigr)   \exp(-\lambda_n)\bigl({\lambda_n^k}/{k!}\bigr)\\
&&\hspace*{88pt}\qquad{} +
\Biggl({ (T+1)^2 }\Big/\Biggl(\prod_{j=1}^T
\bigl(1-({x}/{j})\bigr)\Biggr)\Biggr)\Biggr),\nonumber
\end{eqnarray}
and at most
%
\begin{eqnarray}
\label{usedb2} && \Biggl(\sum_{k=0}^{n-(T+1)}
(n - k) \prod_{j=T+1}^{n-k}\biggl(1-\frac{x}{j}\biggr) \exp(-\lambda_n)
\bigl({\lambda_n^k}/{k!}\bigr)\nonumber\\
&&\hspace*{55pt}\quad{}+ \Biggl({ (T+1)^2 }\Big/{\Biggl(\prod_{j=1}^T \bigl(1-({x}/{j})\bigr)\Biggr)}\Biggr)\Biggr) \\
&&\qquad{}\Big/\Biggl(
\lambda_n^{-{1}/{2}} \sum_{k=0}^{n-(T+1)} k \prod_{j=T+1}^{n-k}\bigl(1-({x}/{j})\bigr)   \exp(-\lambda_n)\bigl({\lambda_n^k}/{k!}\bigr)\Biggr)
.\nonumber
\end{eqnarray}
We now simplify the terms in (\ref{usedb1}) and (\ref{usedb2}), by
proving that
for all $n \geq T+1$, and $k \in[0, n - T - 1]$,
%
\begin{equation}
\label{aslemma1}\qquad \exp\bigl(-2 T^{-1}\bigr) (n-k)^{-x}
T^x \leq\prod_{j=T+1}^{n-k}
\biggl(1 - \frac
{x}{j}\biggr) \leq\exp\bigl(2 T^{-1}\bigr)
(n-k)^{-x} T^x.
\end{equation}
Indeed, since $0 < x < 1$, it follows from a simple Taylor series
expansion that for all $j \geq3$, $1 \leq
\frac{ \exp( - {x}/{j}) }{ 1 - ({x}/{j}) } \leq1 + j^{-2}$.
Thus for $j \geq T+1$,
\[
\prod_{j=T+1}^{n-k} \frac{ \exp(-{x}/{j}) } { (1 -
({x}/{j})) } \leq
\prod_{j = T+1}^{n-k} \bigl(1 + j^{-2}
\bigr) \leq \exp\biggl(\int_{T}^{\infty}
x^{-2} \,dx \biggr) = \exp\bigl(T^{-1}\bigr)
\]
and
%
\begin{equation}
\label{intermediate1}\quad  \exp\bigl(-T^{-1}\bigr) \prod
_{j=T+1}^{n-k}\exp\biggl(-\frac{x}{j}\biggr) \leq
\prod_{j=T+1}^{n-k} \biggl(1 -
\frac{x}{j}\biggr) \leq \prod_{j=T+1}^{n-k}
\exp\biggl(-\frac{x}{j}\biggr).
\end{equation}
Let $H_k \stackrel{\Delta}{=} \sum_{j=1}^k \frac{1}{j}$ denote the
$k$th harmonic number. Then it follows from the results of~\cite{Y91},
and the fact that $n - k > T$, that
%
\begin{equation}
\label{logharm} \log\biggl( \frac{n-k}{T} \biggr) - (2 T)^{-1}
\leq H_{n-k} - H_T \leq\log\biggl( \frac
{n-k}{T} \biggr)
+ (2 T)^{-1}.
\end{equation}
Combining (\ref{intermediate1}) and (\ref{logharm}) with the fact that
$0 < \frac{x}{2T} < (2T)^{-1}$ completes the proof of (\ref{aslemma1}).

It follows from (\ref{usedb1}), (\ref{usedb2}) and (\ref
{aslemma1}) that for all sufficiently large $n$, $ ( z_n(x) -
\lambda_n  )\lambda_n^{-{1}/{2}}$ is at least
%
\begin{eqnarray}
\label{usedc1}&& \exp\bigl(- 4 T^{-1} \bigr) \Biggl( \sum_{k=0}^{n-(T+1)} (n - k)^{1-x} \exp
(-\lambda_n)({\lambda_n^k}/{k!}) \Biggr)\nonumber\\
&&\qquad{} \Big/\Biggl(
\lambda_n^{- {1}/{2}} \sum_{k=0}^{n-(T+1)} k (n-k)^{-x}  \exp
(-\lambda_n)\bigl({\lambda_n^k}/{k!}\bigr) \\
&&\hspace*{41pt}\qquad{}+ \Biggl({(T+1)^2 T^{-x} }\Big/{\Biggl(\prod_{j=1}^T\biggl (1-\frac{x}{j}\biggr)\Biggr)}\Biggr)
\Biggr),\nonumber
\end{eqnarray}
and at most
%
\begin{eqnarray}
\label{usedc2} &&\exp\bigl(4 T^{-1}\bigr) \Biggl( \sum_{k=0}^{n-(T+1)} (n - k)^{1-x} \exp(-\lambda_n)\bigl({\lambda_n^k}/{k!}\bigr)\nonumber\\
&&\hspace*{31pt}\qquad{} +
\Biggl({\bigl((T+1)^2 T^{-x} \bigr)}
\Big/{\Biggl(\prod_{j=1}^T \biggl(1-\frac{x}{j}\biggr)\Biggr)}\Biggr)\Biggr)\\
&&\qquad{} \Big/ \Biggl(
\lambda_n^{-{1}/{2}} \sum_{k=0}^{n-(T+1)} k (n-k)^{-x}   \exp
(-\lambda_n)\bigl({\lambda_n^k}/{k!}\bigr) \Biggr).\nonumber
\end{eqnarray}
With inequalities (\ref{usedc1}) and (\ref{usedc2}) in hand, we are now
in a position to complete the proof of Lemma~\ref{useddd1}.
We begin by proving the lower bound. The term $\lambda_n^{-{1}/{2}}
\sum_{k=0}^{n-(T+1)} k (n-k)^{-x}   \exp(-\lambda_n)\frac{\lambda_n^k}{k!}$ appearing in the denominator of (\ref{usedc1}) is at most
%
\begin{eqnarray}
\label{anothertobound}\quad && \lambda_n^{-{1}/{2}} \sum
_{k=0}^{ \lceil n - T^{-1} n^{{1}/{2}}
\rceil} k (n-k)^{-x} \exp(-
\lambda_n)\frac{\lambda_n^k}{k!} + \lambda_n^{-{1}/{2}}
\max_{0 \leq k \leq n } \biggl( k \exp (-\lambda_n) \frac{\lambda_n^k}{k!}
\biggr)
\nonumber
\\[-8pt]
\\[-8pt]
\nonumber
&&\qquad{}\times \sum_{k = \lceil n - T^{-1} n^{{1}/{2}} \rceil+ 1}^{n - (T+1)} (n-k)^{-x}.
\end{eqnarray}
Recall that for all sufficiently large $n$, $\sup_{k \geq0}  (
\exp
(-\lambda_n) \frac{\lambda_n^k}{k!}  ) \leq n^{-{1}/{2}}$, and
$(\frac{n}{\lambda_n})^{{1}/{2}} \leq2$, from which it follows
that the second summand of (\ref{anothertobound}) is at most
\begin{eqnarray*}
\biggl(\frac{n}{\lambda_n}\biggr)^{{1}/{2}} \sum
_{k = \lceil n - T^{-1} n^{{1}/{2}} \rceil+ 1}^{n - (T+1)} (n-k)^{-x} &\leq& 2 \int
_0^{ T^{-1} n^{{1}/{2}} } y^{-x} \,dy\\
& =& 2
(1-x)^{-1} T^{-(1-x)} n^{{(1-x)}/{2}}.
\end{eqnarray*}
Using the above to upper-bound the denominator of (\ref{usedc1}),
multiplying through by $\lambda_n^{{(x-1)}/{2}}$ and observing that
$\lambda_n^{{(x-1)}/{2}} \frac{(T+1)^2 T^{-x} }{\prod_{j=1}^T
(1-\frac
{x}{j})} \leq2 (1-x)^{-1} T^{-(1-x)}$ for all sufficiently large $n$
completes the proof of the lower bound. The upper bound follows from a
similar argument, and we omit the details.
\end{pf*}

\begin{pf*}{Proof of Lemma~\ref{parabolicconcave}}
We begin by demonstrating that $z_{n,k}$ is a twice-differentiable
concave function on $(0,1)$ for all $k \leq n$, which will imply that
$z_{\infty}$, and ultimately $\upsilon(x,-B)$, are concave by taking
limits. We proceed by induction on $k$. The base case $k = 1$ is
trivial, since $z_{n,1}(x) = \lambda_n + 1 - x$.
Now, let us assume the statement is true for $j = 1,\ldots,k-1$ with
$k-1 \leq n-1$. It may be easily verified that
$f_{n,k}(x) = (\lambda_n + k - x)f_{n,k-1}(x) - \lambda_n(k-1)
f_{n,k-2}(x)$. Thus since $z_{n,k-1}$ is strictly positive on $(0,1)$,
which follows from
Lemma~\ref{kijima2}(ii), we find that
\begin{eqnarray*}
&&\hspace*{-4pt}\frac{d^2}{d x^2} z_{n,k}(x) \\
&&\hspace*{-6pt}\quad= \lambda_n (k-1) \biggl( -2
z_{n,k-1}(x)^{-3} \biggl( \frac{d}{dx} z_{n,k-1}(x)
\biggr)^2 + z_{n,k-1}(x)^{-2} \frac{d^2}{dx^2}
z_{n,k-1}(x) \biggr).
\end{eqnarray*}
Since the induction hypothesis implies that $\frac{d^2}{dx^2}
z_{n,k-1}(x) \leq0$, it follows that $z_{n,k}$ is twice-differentiable
on $(0,1)$ and satisfies $\frac{d^2}{d x^2} z_{n,k}(x) \leq0$
(concavity), proving the induction.

Combining the above with Proposition~\ref{astep1}, and the
fact that pointwise limits of concave functions are concave,
demonstrates that $z_{\infty}$ is a concave function of $x$ on $(0,1)$
for any fixed $B > 0$. Observing that $\upsilon(x,-B) = z_{\infty}(x) -
B$ completes the proof.
\end{pf*}

\begin{pf*}{Proof of Lemma~\ref{interintegralbounds}}
Let $d_n(a) \stackrel{\Delta}{=} 6  ( 1 + s_n(a)  )  (a +
n^{-{1}/{2}} )$. Then the left-hand side of (\ref
{anotherinteqq}) equals
%
\begin{eqnarray}\label{line1}\quad
& & \int_{ (n^{{1}/{2}} - \lambda
_n^{{1}/{2}})^2}^{ 2 ( \lambda_n n )^{{1}/{2}} }
\exp \biggl(\frac{a}{3} d_n(a) - x t + d_n(a)
x^{{1}/{2}} \biggr) \bigl(x - \bigl(n^{{1}/{2}} - \lambda_n^{{1}/{2}}
\bigr)^2 \bigr)^{-{1}/{2}}
\nonumber
\\[-8pt]
\\[-8pt]
\nonumber
&&\hspace*{50pt}{}\times \bigl(\bigl(n^{{1}/{2}} +
\lambda_n^{{1}/{2}}\bigr)^2 - x \bigr)^{-{1}/{2}}
\,dx
\\
\label{line2}&&\qquad{}+  \int_{ 2 (\lambda_n
n)^{{1}/{2}} }^{ (n^{{1}/{2}} + \lambda_n^{{1}/{2}})^2 }
\exp \biggl(\frac{a}{3} d_n(a) - x t + d_n(a)
x^{{1}/{2}} \biggr) \nonumber\\
&&\hspace*{89pt}{}\times\bigl(x - \bigl(n^{{1}/{2}} - \lambda_n^{{1}/{2}}
\bigr)^2 \bigr)^{-{1}/{2}}
\\
&&\hspace*{89pt}{}\times \bigl(\bigl(n^{{1}/{2}} +
\lambda_n^{{1}/{2}}\bigr)^2 - x \bigr)^{-{1}/{2}}
\,dx.\nonumber
\end{eqnarray}
Let $u_n \stackrel{\Delta}{=} 2 (\lambda_n n)^{{1}/{2}} -
(n^{{1}/{2}}-\lambda_n^{{1}/{2}})^2$.
Since $ ((n^{{1}/{2}} + \lambda_n^{{1}/{2}})^2 - x
)^{-{1}/{2}} \leq ( (n^{{1}/{2}} + \lambda_n^{{1}/{2}})^2 - 2 (\lambda_n n)^{{1}/{2}}  )^{-{1}/{2}}$
for $x
\in ( (n^{{1}/{2}} - \lambda_n^{{1}/{2}})^2, 2 (\lambda_n
n)^{{1}/{2}}  )$, (\ref{line1}) is at most
\begin{eqnarray*}
& &\hspace*{-4pt} \exp \biggl(\frac{a}{3} d_n(a) \biggr) \bigl(
\bigl(n^{{1}/{2}} + \lambda_n^{{1}/{2}}\bigr)^2
- 2 (\lambda_n n)^{{1}/{2}} \bigr)^{-{1}/{2}}\\
&&\hspace*{-6pt}\qquad{} \times\int
_{(n^{{1}/{2}} - \lambda_n^{{1}/{2}})^2}^{ 2
(\lambda_n n)^{{1}/{2}} }\exp \bigl(- x t + d_n(a)
x^{{1}/{2}} \bigr) \bigl(x - \bigl(n^{{1}/{2}} - \lambda_n^{{1}/{2}}
\bigr)^2 \bigr)^{-{1}/{2}}\,dx
\\
& &\hspace*{-6pt}\quad = \exp \biggl(\frac{a}{3} d_n(a) \biggr) (
\lambda_n + n )^{-{1}/{2}} \\
&&\hspace*{-6pt}\qquad{}\times\int_{0}^{ u_n }
\exp \bigl(- \bigl(y + \bigl(n^{{1}/{2}}-\lambda_n^{{1}/{2}}
\bigr)^2 \bigr) t + d_n(a) \bigl(y + \bigl(n^{{1}/{2}}-
\lambda_n^{{1}/{2}}\bigr)^2 \bigr)^{{1}/{2}}
\bigr)y^{-{1}/{2}}\,dy
\\
& &\hspace*{-6pt}\quad \leq ( \lambda_n + n )^{-{1}/{2}} \exp \biggl(
\frac
{a}{3} d_n(a) + d_n(a) \bigl(n^{{1}/{2}}-
\lambda_n^{{1}/{2}}\bigr) - \bigl(n^{{1}/{2}}-
\lambda_n^{{1}/{2}}\bigr)^2 t \biggr)\\
&&\hspace*{-6pt}\hspace*{20pt}{}\times \int
_{0}^{ u_n } \exp \bigl(- y t + d_n(a)
y^{{1}/{2}} \bigr)y^{-{1}/{2}}\,dy,
\end{eqnarray*}
where the final inequality follows from the fact that $ (y +
(n^{{1}/{2}}-\lambda_n^{{1}/{2}})^2 )^{{1}/{2}}
\leq
y^{{1}/{2}} + n^{{1}/{2}}-\lambda_n^{{1}/{2}}$. It may be
easily verified that
$- y t + d_n(a) y^{{1}/{2}} \leq- \frac{1}{2} y t +
d_n(a)^2(2t)^{-1}$ for all $y > 0$, and $\int_{0}^{ \infty} \exp(-
\frac{1}{2} y t)y^{-{1}/{2}}\,dy = (\frac{2 \pi}{t})^{{1}/{2}}$,
and we conclude that (\ref{line1}) is at most
%
\begin{eqnarray}
\label{firsttermbound} &&J_1 \stackrel{\Delta} {=} \biggl(
\frac{\pi}{\lambda_n t}\biggr)^{{1}/{2}}  \exp \biggl( \frac{a}{3}
d_n(a) + d_n(a) \bigl(n^{{1}/{2}}-
\lambda_n^{{1}/{2}}\bigr) + d_n(a)^2(2t)^{-1}
\nonumber
\\[-8pt]
\\[-8pt]
\nonumber
&&\hspace*{211pt}{}- \bigl(n^{{1}/{2}}-\lambda_n^{{1}/{2}}\bigr)^2
t \biggr).
\end{eqnarray}
We now bound (\ref{line2}). Let $S \stackrel{\Delta}{=}  ( 2 (
\lambda_n n )^{{1}/{2}} , (n^{{1}/{2}} + \lambda_n^{{1}/{2}})^2  )$.
Since $x \in S$ implies
\[
\bigl(x - \bigl(n^{{1}/{2}} - \lambda_n^{{1}/{2}}
\bigr)^2 \bigr)^{-{1}/{2}} \leq \bigl( 2 (\lambda_n
n)^{{1}/{2}} - \bigl(n^{{1}/{2}} - \lambda_n^{{1}/{2}}
\bigr)^2 \bigr)^{-{1}/{2}} \leq ( 3 \lambda_n - n
)^{-{1}/{2}},
\]
(\ref{line2}) is at most
\begin{eqnarray*}
& & \exp \biggl( \frac{a}{3} d_n(a) \biggr) ( 3
\lambda_n - n )^{-{1}/{2}} \sup_{z \in S} \exp \bigl(- z t +
d_n(a) z^{{1}/{2}} \bigr)\\
&&\quad{}\times \int_{2 (\lambda_n n)^{{1}/{2}} }^{(n^{{1}/{2}} + \lambda
_n^{{1}/{2}})^2}
\bigl(\bigl(n^{{1}/{2}} + \lambda_n^{{1}/{2}}
\bigr)^2 - x \bigr)^{-{1}/{2}}\,dx
\\
& &\qquad = \exp \biggl( \frac{a}{3} d_n(a) \biggr) ( 3
\lambda_n - n )^{-{1}/{2}} \sup_{z \in S} \exp \bigl(- z t +
d_n(a) z^{{1}/{2}} \bigr)\\
&&\quad\qquad{}\times \int_{0}^{\lambda_n + n}
y^{-{1}/{2}} \,dy,
\end{eqnarray*}
which is itself at most
%
\begin{equation}
\label{anotherinter111} J_2 \stackrel{\Delta} {=} 2\biggl(
\frac{\lambda_n + n}{3 \lambda_n -
n}\biggr)^{{1}/{2}} \exp \biggl( \frac{a}{3}
d_n(a) + d_n(a)^2(2t)^{-1} - (
\lambda_n n)^{{1}/{2}} t \biggr),
\end{equation}
where the final inequality follows from the fact that $- z t + d_n(a)
z^{{1}/{2}} \leq-\frac{1}{2} z t + d_n(a)^2(2t)^{-1}$, and $\int_{0}^{\lambda_n + n} y^{-{1}/{2}} \,dy = 2 (\lambda_n + n)^{{1}/{2}}$.
It may be easily verified that there exists $N_{B,a}, C_{B,a} < \infty
$, depending only on $B$ and $a$, s.t. for all $n \geq N_{B,a}$ and $t
\geq1$,
one has $d_n(a) \leq6 a + C_{B,a} n^{-{1}/{2}}$, $n^{{1}/{2}}
- \lambda_n^{{1}/{2}} \leq\frac{B}{2} + C_{B,a} n^{-{1}/{2}}$ and
$J_2 \leq n^{-1} J_1$.
The lemma then follows by using (\ref{firsttermbound}) to bound~(\ref
{line1}), (\ref{anotherinter111}) to bound~(\ref{line2}) and applying a
simple Taylor series expansion.
\end{pf*}
\end{appendix}0

\section*{Acknowledgments}
The authors would like to thank Hans Blanc, Ton Dieker, Erik van Doorn,
Peter Glynn, Johan van Leeuwaarden, Bill Massey, Josh Reed and Ward Whitt
for their helpful discussions and insights. The authors especially
thank Johan van Leeuwaarden for his insights into the parabolic
cylinder functions, and sharing an early draft of his work. The authors
also thank two anonymous referees, who helped to significantly improve
the paper's presentation.

%
%


\printaddresses

\end{document}